\documentclass[12pt]{article}

\usepackage{amssymb, amsmath}
\usepackage{ulem}
\usepackage{amsfonts}
\usepackage{graphicx}
\usepackage{setspace} 
\usepackage{rotating}
\usepackage{enumitem}
\allowdisplaybreaks
\usepackage{hyperref}

\usepackage{tikz}
\usepackage{mathrsfs}
\usepackage{url}
\usepackage[authoryear]{natbib}
\newcommand{\be}{\begin{equation}}
	\newcommand{\ee}{\end{equation}}
\newcommand{\bea}{\begin{eqnarray}}
	\newcommand{\eea}{\end{eqnarray}}
\newcommand{\bean}{\begin{eqnarray*}}
	\newcommand{\eean}{\end{eqnarray*}}
\newcommand{\brray}{\begin{array}}
	\newcommand{\erray}{\end{array}}
\newcommand{\ben}{\begin{equation}{nonumber}}
	\newcommand{\een}{\end{equation}{nonumber}}

\newtheorem{dfn}{Definition}[section]
\newtheorem{thm}[dfn]{Theorem}
\newtheorem{lmma}[dfn]{Lemma}
\newtheorem{ppsn}[dfn]{Proposition}
\newtheorem{crlre}[dfn]{Corollary}
\newtheorem{xmpl}[dfn]{Example}
\newtheorem{rmrk}[dfn]{Remark}

\newcommand{\bdfn}{\begin{dfn}}
	\newcommand{\bthm}{\begin{thm}}
		\newcommand{\blmma}{\begin{lmma}}
			\newcommand{\bppsn}{\begin{ppsn}}
				\newcommand{\bcrlre}{\begin{crlre}}
					\newcommand{\bxmpl}{\begin{xmpl}}
						\newcommand{\brmrk}{\begin{rmrk}}
							\newcommand{\edfn}{\end{dfn}}
						\newcommand{\ethm}{\end{thm}}
					\newcommand{\elmma}{\end{lmma}}
				\newcommand{\eppsn}{\end{ppsn}}
			\newcommand{\ecrlre}{\end{crlre}}
		\newcommand{\exmpl}{\end{xmpl}}
	\newcommand{\ermrk}{\end{rmrk}}

\newcommand{\clb}{{\cal B}}

\newtheorem{theorem}{Theorem}[section]
\newtheorem{example}{Example}[section]
\newtheorem{remark}{Remark}[section]
\newtheorem{lemma}{Lemma}[section]
\newtheorem{definition}{Definition}[section]

\newtheorem{proposition}{Proposition}[section]

\def\a*{{\cal A}_{h,*}}
\def\B{{\cal B}(h)}
\def\B1{{\cal B}_1(h)}
\def\b{{\cal B}^{\rm s.a.}(h)}
\def\b1{{\cal B}^{\rm s.a.}_1(h)}

\usepackage{geometry}
\usepackage{xr}
\makeatletter
\newcommand*{\addFileDependency}[1]{
  \typeout{(#1)}
  \@addtofilelist{#1}
  \IfFileExists{#1}{}{\typeout{No file #1.}}
}
\makeatother

\newcommand*{\myexternaldocument}[1]{
    \externaldocument{#1}
    \addFileDependency{#1.tex}
    \addFileDependency{#1.aux}
}
\myexternaldocument{Main_Bernoulli}
\begin{document}




  \title{Co-variance Operator of Banach Valued Random Elements: U-Statistic Approach}
\author{\small 
	Suprio Bhar \\
	\small IIT Kanpur\\
	\small Department of Mathematics and Statistics \\
	\small  Kanpur 208016, India\\
	{\small email: suprio@iitk.ac.in}\\
	\and
	\small Subhra Sankar Dhar \\
	\small  IIT Kanpur\\
	\small   Department of Mathematics and Statistics \\
	\small Kanpur 208106, India\\
	{\small email: subhra@iitk.ac.in}\\
}
\date{}
\maketitle

\setcounter{equation}{0}

\setcounter{section}{0}\ \\


\begin{abstract}
 This article proposes a co-variance operator for Banach valued random elements using the concept of $U$-statistic. We then study the asymptotic distribution of the proposed co-variance operator along with related large sample properties. Moreover, specifically for Hilbert space valued random elements, the asymptotic distribution of the proposed estimator is derived even for dependent data under some mixing conditions. Finally, a small numerical study gives some ideas about the choice of the tuning parameter involved in the proposed co-variance operator. 
\end{abstract}


\section{Introduction}\label{Intro}
Since the last decade or so, there has been considerable attention on infinite dimensional data in Statistics and related subjects' literature as nowadays in many applications, the dimension of the data is larger than the sample size and exhibits a certain degree of smoothness, which can be embedded into an appropriate infinite dimensional space. To analyze such data, one may consider that the data/observations are realizations of a random element $X$ defined on the Banach space ${\cal{B}}$ or on the certain Hilbert space ${\cal{H}}$, which are all infinite dimensional in nature. For example, functional data (see, e.g., \cite{Ramsey2002} and  \cite{Fer2006}) is such type of data, and to study the functional data, one may adopt the techniques of infinite dimensional space. Though adopting the usual technique of multivariate analysis may be possible by observing the functions on discrete time points, it may fail to capture the smoothness of the functions. Moreover, the characteristic parameters (e.g., mean/Expectation) of the measure associated with the random element are generally infinite dimensional themselves, and hence, one needs to consider the mathematical technique used in infinite dimensional space. 

In this article, we investigate the second order properties of the infinite dimensional random element, and as in the case of finite dimensional random element, covariance {\it operator} is the second order characterization of a random element defined on ${\cal B}$. For literature survey on the covariance operator, the readers may refer to \cite{MR0517336}, \cite{MR3127875}, \cite{MR3127883}, \cite{MR3556768}, \cite{MR3771767}, \cite{MR3809475}, \cite{MR3984017}, \cite{MR4036049}, \cite{MR4259441}, \cite{MR4509065}, \cite{MR4499013}, \cite{MR4361616}, \cite{MR4411501}, \cite{MR4517351}, \cite{HaoyuWang} and a few references therein. In view of applications also, the inference on covariance operator is an integral part of many cases. For example, it may be useful to detect whether the data has any outliers or not OR it can be used in functional principal component analysis, which is a common tool in dimension reduction. In the context of a real life example, \cite{Panaretos2010} studied DNA minicircles and found no differences in their means whereas two different groups have different covariance structures. Precisely speaking, the covariance structure/operator is associated to the flexibility or stiffness of the DNA minicircles, and two treatment groups have differences in terms of flexibility or stiffness. 
For some other real life examples associated with covariance operators, the readers may refer to \cite{10.1214/17-EJS1347}, \cite{PMID:18580414}, \cite{COFFEY20111144} and a few relevant references therein. 

There have been a few attempts related to statistical inference on covariance operators. For instance, \cite{FERRATY20074903} investigated the features of several curves based on comparing covariance operators. It follows from Kosambi-Karhunen–Loève expansion (see, e.g., \cite{Kosambi1943}, \cite{MR0023013} and \cite{MR0651018}) that their proposed methodology is equivalent to testing whether all the samples have the same set of functional principal components of the covariance operators or not. Specifically, for two-sample hypothesis problems, \cite{10.1214/07-AOS516}, \cite{Fremdt2013} and \cite{Panaretos2010} used the similar idea to formulate the test statistic. Besides, \cite{RePEc:spr:aistmt:v:74:y:2022:i:2:d:10.1007_s10463-021-00795-2} proposed statistical inference tools for the covariance operators of functional time series in the two sample and change point problem. Their approach is not testing the null hypothesis of exact equality of the covariance operators. Instead, they proposed to formulate the null hypotheses in the form that ``the distance between the operators is small", where they measure deviations by the sup-norm.  

Moreover, there are a few articles on robustification of covariance operators as well. \cite{10.1093/biomet/ass037} introduced the notion of a covariance operator, investigate its use in probing the second-order structure of functional data, and develop a test for comparing the second-order characteristics of two functional samples that is resistant to a typical observations and departures from normality. In this spirit, \cite{BOENTE2019115} studied the asymptotic behavior of the sample spatial sign covariance operator centered at an estimated location. Furthermore,   \cite{RePEc:spr:aistmt:v:70:y:2018:i:4:d:10.1007_s10463-017-0613-1} extended the classical two-population problem, presenting a test for equality of covariance operators among k ($\geq 2$) populations in which the asymptotic distribution of the sample covariance operator plays a crucial role in deriving the asymptotic distribution of the proposed statistic. It is well known that the presence of outliers in the sample might lead to invalid conclusions. All these works motivate us to study a class of covariance operators, which can be robust and efficient as well for various choices of a certain tuning parameter. 

This article proposes a new variant of covariance operator using the concept of $U$-statistic of Banach valued random elements (see \cite{Borovskikh1996}), and different choices of the order of the kernel (denoted by $m$) associated with $U$-statistic provide different variants of covariance operator. To investigate the performance of different variants, in Section \ref{NS}, a small numerical study reveals that the choices of the tuning parameter, i.e., $m$, can affect the efficiency of the covariance operator for various distributions. Moreover,  in the course of this study, we establish results related to the large sample of the proposed co-variance operator. Additionally, for Hilbert space valued {\it dependent} data, the asymptotic distribution of the proposed operator is derived. Finally, it is to be noted that the classical co-variance operator can be obtained when $m = 1$, and hence, the large sample results related to the classical covariance operator can be derived from the results obtained in this article.


In the course of studying all these aforementioned large sample properties of the proposed operator, we face a few mathematical challenges. First and foremost, \cite{Borovskikh1996} defined $U$-statistic for ${\cal{B}}$-valued random elements, where ${\cal{B}}$ is a Banach space. However, our proposed co-variance operator is based on $(\cal B \otimes \cal B)$-valued random element, and hence, one cannot directly use the results of $U$-statistic for ${\cal{B}}$-valued random elements. In this work, this issue is dealt using advanced techniques in functional analysis as $(\cal B \otimes \cal B)$ is not even uniquely defined (see, e.g., \cite{RyanBk}; see also Appendix for details on $(\cal B \otimes \cal B)$). Secondly, as the probability theory of $U$-statistic for ${\cal{B}}$-valued random elements depends on the geometry of the Banach space ${\cal{B}}$ and the rank of the kernel associated with the $U$-statistic, one needs to take into account the geometry of tensor space $(\cal B \otimes \cal B)$ to establish all these results. Thirdly, establishing the similar results for dependent data, one needs to impose a certain structure on the Banach space ${\cal{B}}$ and the mixing conditions on the random elements. In this work, we overcome all these issues and establish all relevant results. 

The lack of an inner-product operation on a Banach space, compared to any Hilbert space, creates an ambiguity while trying to understand `angles'. In the absence of clear geometric structures on a Banach space, many of the standard probabilistic methodologies, developed for random elements taking values in finite-dimensional Euclidean spaces or infinite dimensional separable Hilbert spaces, does not have straight-forward generalizations for random elements taking values in some Banach space. It is well-known that the extensions of the Law of Large Numbers or the Central Limit Theorem for random elements in Banach spaces require the Banach spaces to satisfy certain geometric conditions, typically described in terms of type and co-types (see \cite{WoyczynskiBk}). In the Appendix, we briefly recall some of these concepts on Banach spaces, which we use in this article. The interested reader may also see further references on this topic, such as \cite{RyanBk, Chatterji_1968, Borovskikh1996, Ledoux-TalagrandBk,MetivierBk} and the references therein.

The rest of the article is organized as follows. In Section \ref{PCO}, the new version of co-variance operator using the idea of $U$-statistic is proposed. Section \ref{LSP} studies various large sample properties of the proposed co-variance operator, and a small numerical study is conducted about the choice of the tuning parameter $m$ in Section \ref{NS}. Section \ref{CR} consists of a few concluding remarks, and finally, Section \ref{AP} contains all technical details. 

\section{Proposed Co-variance Operator}\label{PCO}
Let ${\cal{B}}$ be a real separable Banach space with a norm $||.||_{{\cal{B}}}$, and let ${\cal{B}}^{*}$ denote the dual space of ${\cal{B}}$. Suppose that $X_{1}, \ldots, X_{n}$ (identically distributed with $X$) are ${\cal{B}}$-valued i.i.d.\ random elements with probability law $P$, and the corresponding measurable space is $(\cal{X}, {\cal{A}})$. Under this set-up, we propose the following sample covariance operator for any $1\leq m\leq n$. 
\begin{equation}\label{Sampleversion}
C_{m, n} = \frac{1}{{n\choose m}}\sum_{1\leq i_{1}<i_{2}<\ldots<i_{m}\leq n} \left(\frac{X_{i_{1}} + \ldots + X_{i_{m}}}{m} - \theta (m)\right)\otimes\left(\frac{X_{i_{1}} + \ldots + X_{i_{m}}}{m} - \theta (m)\right),  
\end{equation} where 
$$\theta (m) = \theta (P; m) = E\left(\frac{X_{1} + \ldots + X_{m}}{m}\right).$$ Here $\otimes$ is an appropriate notion of tensor product, which may not be unique for Banach valued variables. In this work, we consider two types of tensor products, namely, projective tensor product and injective tensor product (see, e.g., \cite{RyanBk}). All results stated in the subsequent sections are valid for both aforementioned tensor products. In the course of this study, in many places, we put assumptions directly on the tensor space $\mathcal{B}\otimes\mathcal{B}$ rather than on $\mathcal{B}$ as characterization of the geometric feature of $\mathcal{B}\otimes\mathcal{B}$ based on that of $\mathcal{B}$ is not always tractable (see, e.g., \cite{MR4477944}) Here, in the expression of $\theta(m)$, $E(.)$ is taken as Bochner sense (see, e.g., \cite{Bochner1933}), and note that $\theta(m) = \theta(P; m) = E(X_{1})$ if $\{X_{1}, \ldots, X_{n}\}$ are identically distributed random variables. The population version of the corresponding covariance operator is the following. 
\begin{equation}\label{Populationversion}
C_{m} = E\left\{\left(\frac{X_{1} + \ldots + X_{m}}{m} - \theta (m)\right)\otimes\left(\frac{X_{1} + \ldots + X_{m}}{m} - \theta (m)\right)\right\}.   
\end{equation} 
Note that $C_{m, n}$ is a $\mathcal{B}\otimes\mathcal{B}$-valued random elements, and $C_{m}$ is a $\mathcal{B}\otimes\mathcal{B}$-valued variable.  Proposition \ref{propu} asserts how $C_{m, n}$ and $C_m$ are directly associated. 

\begin{ppsn}\label{propu}
$E(C_{m, n}) = C_{m}$ for all $n\in\mathbb{N}$ if $E||X||^{2}_{\cal{B}} < \infty$.
\end{ppsn}

It is shown in Proposition \ref{propu} that $C_{m, n}$ is an unbiased estimator of $C_{m}$ for all $n\in\mathbb{N}$ if $E||X||^{2}_{\cal{B}} < \infty$. Note that though $C_{m, n}$ is a ${\cal{B}}\otimes {\cal{B}}$-valued random element, the required condition on $C_{m, n}$ to be the unbiased estimator of $C_{m}$ involves the moment condition on ${\cal{B}}$-valued random element $X$. Hence, one does not need to check any condition associated with the complicated tensor product to use $C_{m, n}$ as an unbiased estimator of $C_{m}$.

\begin{remark}\label{Cmn-justification}
We now want to discuss why $C_{m, n}$ or $C_{m}$ can be considered as a measure of dispersion. First note that $C_{m}$ coincides with the usual covariance when $m = 1$ and ${\cal{B}} = \mathbb{R}$, which one would hope to see. For general $m$ and ${\cal{B}} = \mathbb{R}$, one can view $C_{m}$ as $$C_{m} = E\left(\frac{X_{1} - \theta(m)}{m} + \ldots + \frac{X_{m} - \theta(m)}{m}\right)^{2},$$ i.e., $C_{m}$ measures a certain dispersion from $X_{i}$ to $\theta(m)$ ($i = 1, \ldots, m$), and hence, one may consider $C_{m}$ (or $C_{m, n}$) as a certain measure of dispersion.
\end{remark}

\begin{remark}\label{Cmn-Hilbert}
In the case of a real separable Hilbert space ${\cal{H}}$, elements of the tensor product are identified as scalar valued bilinear maps on ${\cal{H}} \times {\cal{H}}$, and the inner-product on $\cal{H} \otimes \cal{H}$ is defined using the inner-product $\langle \cdot, \cdot \rangle_{\cal{H}}$ as follows,
\[\langle x_1 \otimes x_2, h_1 \otimes h_2 \rangle_{\cal{H}\otimes \cal{H}} := \langle x_1, h_1 \rangle_{\cal{H}} \langle x_2, h_2 \rangle_{\cal{H}},\]
for all $x_1, x_2, h_1, h_2 \in \cal{H}$. Without loss of generality, taking $\theta(m) = 0$, we now have
\begin{align*}
C_{m, n} &= \frac{1}{{n\choose m}} \frac{1}{m^2}\sum_{1\leq i_{1}<i_{2}<\ldots<i_{m}\leq n} \left(X_{i_{1}} + \ldots + X_{i_{m}}\right)\otimes\left(X_{i_{1}} + \ldots + X_{i_{m}}\right)\\
&= \frac{1}{{n\choose m}} \frac{1}{m^2} \sum_{1\leq i_{1}<i_{2}<\ldots<i_{m}\leq n}\quad \sum_{i, j \in \{i_1, i_2, \cdots, i_m\}} X_{i} \otimes X_{j}, 
\end{align*}
and hence, for all $h_1, h_2 \in \cal{H}$,
\begin{align*}
&\langle C_{m, n}, h_1 \otimes h_2 \rangle_{\cal{H}\otimes \cal{H}}\\
&= \frac{1}{{n\choose m}} \frac{1}{m^2} \sum_{1\leq i_{1}<i_{2}<\ldots<i_{m}\leq n}\quad \sum_{i, j \in \{i_1, i_2, \cdots, i_m\}} \langle X_{i} \otimes X_{j} , h_1 \otimes h_2 \rangle_{\cal{H}\otimes \cal{H}}\\
&= \frac{1}{{n\choose m}} \frac{1}{m^2} \sum_{1\leq i_{1}<i_{2}<\ldots<i_{m}\leq n}\quad \sum_{i, j \in \{i_1, i_2, \cdots, i_m\}} \langle X_{i}, h_1 \rangle_{\cal{H}} \langle X_{j}, h_2 \rangle_{\cal{H}}\\
&= \frac{1}{{n\choose m}} \sum_{1\leq i_{1}<i_{2}<\ldots<i_{m}\leq n} \left\{ \frac{1}{m} \sum_{i \in \{i_1, i_2, \cdots, i_m\}} \langle X_{i}, h_1 \rangle_{\cal{H}}  \right\} \left\{ \frac{1}{m} \sum_{j \in \{i_1, i_2, \cdots, i_m\}} \langle X_{j}, h_2 \rangle_{\cal{H}}  \right\}.
\end{align*}
The above action of $C_{m, n}$ on arbitrary $h_1 \otimes h_2$ suggests $C_{m, n}$ as a candidate for the relevant covariance operator in the Hilbert space setting.
\end{remark}

\begin{remark}\label{operator}
 Let us now look at the formulation of $C_{m}$ from the point of view of functional analysis. For any $x\in {\cal{B}}$ and $u\in {\cal{B}}^{*}$, let $u(x)$ denote the action of the linear functional $u$ on $x$. Suppose that $X_{1}, \ldots, X_{n}$ (identically distributed with $X$) are ${\cal{B}}$-valued i.i.d.\  random elements with common probability law $P$. Then, $C_m$ is identified as a real valued bounded bilinear form on $\cal{B}^\ast \times \cal{B}^\ast$ as
 \[C_m(u, v) = E u\left(\frac{X_{1} + \ldots + X_{m}}{m} - \theta (m)\right) v\left(\frac{X_{1} + \ldots + X_{m}}{m} - \theta (m)\right),\]
 for all $u, v\in{\cal{B}}^\ast$. In the same spirit, the sample version of the covariance operator $C_{m, n}$ has the following identification
 \[C_{m, n} (u, v) =  \frac{1}{{n\choose m}}\sum_{1\leq i_{1}<i_{2}<\ldots<i_{m}\leq n} u\left(\frac{X_{i_{1}} + \ldots + X_{i_{m}}}{m} - \theta (m)\right) v\left(\frac{X_{i_{1}} + \ldots + X_{i_{m}}}{m} - \theta (m)\right),\]
for all $u, v\in{\cal{B}}^\ast$.

\end{remark}

In the next section, various large sample statistical properties of this co-variance operator is studied. 

\section{Large Sample Properties}\label{LSP}
This section studies various large sample properties of $C_{m, n}$. Theorem \ref{L1} describes $L_1$ norm convergence of $C_{m, n}$ to $C_{m}$ as $n\rightarrow\infty$, and this implies that $C_{m, n}$ converges to $C_{m}$ in probability as $n\rightarrow\infty$ (see Corollary \ref{Consistency}).   

\begin{theorem}\label{L1}
If $E||X||^{2}_{\cal{B}} < \infty$, then $E||C_{m, n} - C_{m}||_{\cal{B}\otimes {\cal{B}}}\rightarrow 0$ as $n\rightarrow\infty$.
\end{theorem}

\begin{crlre}\label{Consistency}
If $E||X||^{2}_{\cal{B}} < \infty$, then $||C_{m, n} - C_{m}||_{\cal{B}\otimes{\cal{B}}}\stackrel{p}\rightarrow 0$ as $n\rightarrow\infty$.
\end{crlre}

In statistical point of view, the assertions in Theorem \ref{L1} and Corollary \ref{Consistency} along with Proposition \ref{propu} indicate that $C_{m, n}$ can be considered as a reasonably good estimator of $C_{m}$. Now, one may be interested to know the rate of convergence associated with the limits of $C_{m, n}$. Theorem \ref{Rate_of_Convergence} describes the rate of convergence of $C_{m, n}$ to $C_{m}$. Let us first consider the following conditions. 

\vspace{0.1in}

\noindent (C1) $X_{1}, \ldots, X_{n}$ are i.i.d.\ seperable ${\cal{B}}$-valued random elements, where ${\cal{B}}$ is such that 
${\cal{B}}\otimes{\cal{B}}$ is a Banach space of type $p$ for some $p\in (1, 2]$. 

\vspace{0.1in}

\noindent (C2)  For any $m\in\{1, \ldots, n\}$ and any $r\in\{1, \ldots, m\}$, $$E(||X_1 + \ldots + X_{m} - m \theta (m)||^{2}_{\cal{B}} | X_{1}, \ldots, X_{r - 1}) = 0$$ and 
$$Var (E(||X_1 + \ldots + X_{m} - m \theta (m)||^{2}_{\cal{B}} | X_{1}, \ldots, X_{r})) > 0.$$

\vspace{0.1in}

\noindent (C3) For all $q\in (1, 2)$ ($q < p$) and $r\in\{1, \ldots, m\}$, $E||X_1||_{\cal{B}}^{\frac{2mq}{q(m - r) + r}} < \infty$.

\begin{theorem}\label{Rate_of_Convergence}
 Under (C1), (C2) and (C3), $${n\choose r}^{1 - \frac{1}{q}}(C_{m, n} - C_{m})\rightarrow 0$$ almost surely, as $n\rightarrow\infty$ in ${\cal{B}}\otimes{\cal{B}}$. 
\end{theorem}

The interpretation of the conditions (C1), (C2) and (C3) is explained in Remark \ref{discussion_theorem}. As said before, Theorem \ref{Rate_of_Convergence} asserts the rate of convergence of $C_{m, n}$ to $C_{m}$ is ${n\choose r}^{1 - \frac{1}{q}}$, which indicates that for $r = n - 1$ or 1 and $q = 2 - \delta$ for some $\delta > 0$, the rate of convergence will be $n^{\frac{1 - \delta}{2 - \delta}}$, and it will tend to $\sqrt{n}$ if $\delta\rightarrow 0+$. This further implies that the rate of convergence of $(C_{m, n} - C_{m})$ in ${\cal{B}\otimes{\cal{B}}}$ can be made arbitrary close to $\sqrt{n}$. 

\subsection{Asymptotic distribution : Non-degenerate Case}\label{ADNG}

Theorem \ref{asymptotic_normality}states the asymptotic distribution of $(C_{m, n} - C_{m})$ after appropriate normalization with the following conditions. 

\vspace{0.1in}

\noindent (C1*) $X_{1}, \ldots, X_{n}$ are i.i.d.\  seperable ${\cal{B}}$-valued random elements, where ${\cal{B}}$ is such that  $\cal B\otimes \cal B$ is a Banach space of type 2 and $p$-uniformly smooth for $\frac{4}{3}\leq p\leq 2$.

\vspace{0.1in}

\noindent (C2*) 
$Var (E(||X_1 + \ldots + X_{m} - m \theta (m)||^{2}_{\cal{B}} | X_{1})) > 0$.

\vspace{0.1in}

\noindent (C3*) $E||X_1||_{{\cal{B}}}^{\frac{2r}{2r - 1}} < \infty$ for some $r\in\{1, \ldots, m\}$.

\begin{theorem}\label{asymptotic_normality}
Under (C1*), (C2*) and (C3*), $\frac{\sqrt{n}}{m} (C_{m, n} - C_{m})$ converges weakly to $\cal B\otimes \cal B$-valued Gaussian random element $\tau$ having the characteristic function 
$$E(e^{i x^{*}(\tau)}) = e^{-\frac{1}{2} S(x^{*}, x^{*})}, $$ where $x^{*}\in{\cal{B}}^{*}$, and for any $x^{*}\in{\cal{B}}^{*}$ and $y^{*}\in{\cal{B}}^{*}$, $$S(x^{*}, y^{*}) = \int\limits_{{\cal{X}}} x^{*}(g_{1}(z)) y^{*}(g_{1}(z))P(dz).$$ Here $$g_{1}(X_{1}) = E\left[\left\{\left(\frac{X_{1} + \ldots + \ldots + \ldots + X_{m}}{m} - \theta (m)\right)\otimes\left(\frac{X_{1} + \ldots + \ldots + \ldots + X_{m}}{m} - \theta (m)\right)\right\}\vert X_{1}\right].$$ 
\end{theorem}

The assertion in Theorem \ref{asymptotic_normality} indicates that $\frac{\sqrt{n}}{m} (C_{m, n} - C_{m})$ converges weakly to $\cal B\otimes \cal B$-valued Gaussian random element $\tau$ such that for any $x^{*}\in\mathcal{B}^{*}$, $x^{*}(\tau)$ follows a univariate normal distribution with mean $= 0$ and variance $= S(x^{*}, x^{*})$. Observe that this result will enable us to derive the asymptotic efficiency of $C_{m, n}$ for various choices of $m$. Technically speaking, suppose that $S_{1}(x^{*}, y^{*})$ and $S_{2}(x^{*}, y^{*})$ are the asymptotic covariance operators of $C_{m_{1}, n}$ and $C_{m_{2}, n}$ after appropriate normalization for any $x^{*}\in\mathcal{B}^{*}$ and $y^{*}\in\mathcal{B}^{*}$, and in that case, the asymptotic efficiency of $C_{m_{1}, n}$ relative to $C_{m_{2}, n}$ can be defined as $\frac{S_{2}(x^{*}, y^{*})}{S_{1}(x^{*}, y^{*})}$. 

\begin{remark}\label{discussion_theorem}
Here we discuss various conditions assumed in Theorems \ref{Rate_of_Convergence} and \ref{asymptotic_normality}. The condition on the geometry of the space described in (C1) of Theorem \ref{Rate_of_Convergence} is applicable for many well-known spaces, though geometry of tensor products of two or more Banach space is complicated. However, when $\mathcal{B} = \mathcal{H}$, where $\mathcal{H}$ is a some Hilbert space, the tensor space ${\cal{H}}\otimes\mathcal{H}$ will also be type 2 space, and hence, (C1) is applicable on well-known infinite dimensional space like $L_{2}[0, 1]$, $l_{2}[0, 1]$ and many others. For the similar reason, those spaces satisfy the condition (C1*) in Theorem \ref{asymptotic_normality}. The assumptions (C2) and (C2*) in Theorems \ref{Rate_of_Convergence} and \ref{asymptotic_normality}, respectively explains the order of degeneracy of the kernel involved in $C_{m}$, i.e., $\left(\frac{X_{1} + \ldots + X_{m}}{m} - \theta(m)\right)\otimes\left(\frac{X_{1} + \ldots + X_{m}}{m} - \theta(m)\right)$. In order to derive the optimum rate of convergence of $(C_{m, n} - C_{m})$, (C2) precisely indicate that one needs to assume the order of degeneracy of $\left(\frac{X_{1} + \ldots + X_{m}}{m} - \theta(m)\right)\otimes\left(\frac{X_{1} + \ldots + X_{m}}{m} - \theta(m)\right)$ equals with $(r - 1)$, and this $r (\geq 1)$ further involves in the moment assumption in (C3). For instance, let $m = 2$, $r = 1$ and $q = \frac{3}{2}$, (C3) indicates that one needs to have $E||X_{1}||_{\mathcal{B}}^{\frac{12}{5}} < \infty$. To summerize, Theorem \ref{Rate_of_Convergence} asserts that the rate of degeneracy of the Kernel of $C_{m}$, i.e., $r$ controls required moment assumption desribed in (C3) and the rate of convergence of $(C_{m, n} - C_{m})$. The conditions (C2*) and (C3*) of Theorem \ref{asymptotic_normality} implies that one needs to assume that the kernel of $C_{m}$ is non-degenerate,  and the fourth or the lower order moment of the normed random variable is finite. These two assumptions are satisfied for many probability laws defined in infinite dimensional spaces including Euclidean space. 
\end{remark}  


\subsection{Asymptotic distribution : Degenerate Case}\label{DGC}
As it is mentioned in Remark \ref{discussion_theorem} that one needs to assume non-degeneracy of the kernel of $C_{m}$ to have the asymptotic normality of $C_{m, n}$ after a certain normalization. Now, one may be interested to know the asymptotic distribution of $C_{m, n}$ when the kernel of $C_{m}$ is degenerate of a certain order. Theorem \ref{degenerateasymptoic} explores this issue, and for sake of understanding this result, let us define the following.
\begin{equation}\label{Ito}
J_{r}(g_{r}) = \int\limits_{{\cal{X}}}\ldots\int\limits_{{\cal{X}}} g_{r} (x_1, \ldots, x_r) W(dx_1)\ldots W(dx_r),
\end{equation} where
\begin{eqnarray*} &&g_{r} (x_1, \ldots, x_r)\\
&=&  E\left[\left\{\left(\frac{X_{1} + \ldots  + X_{m}}{m} - \theta (m)\right)\otimes\left(\frac{X_{1} + \ldots + X_{m}}{m} - \theta (m)\right)\right\}|X_{1}, \ldots, X_{r}\right],
\end{eqnarray*} and $W$ is a Gaussian random measure on the measurable space $({\cal{X}}, {\cal{A}})$ with mean value zero and for any $A\in{\cal{A}}$ and $B\in{\cal{A}}$, the covariance is given by $E W(A) W(B) = P(A\cap B)$. 


\begin{theorem}\label{degenerateasymptoic} Under (C1), (C2) and (C3*), $n^{\frac{r}{2}}(C_{m, n} - C_{m})$ converges weakly to ${m\choose r}J_{r}(g_{r})$ as $n\rightarrow\infty$. Here $J_{r}(g_{r})$ is the same as \eqref{Ito}.
\end{theorem}

It follows from Theorem \ref{degenerateasymptoic} that the rate of convergence of $(C_{m, n} - C_{m})$, which equals with $n^{\frac{r}{2}}$, depends on the order of degeneracy, i.e., $r$, of the kernel of $C_{m}$. As $r$ increases, the rate of convergence $n^{\frac{r}{2}}$ becomes faster, and consequently, the fastest rate of convergence can be $n^{\frac{m}{2}}$ when $r$
equals with the largest possible value $= m$.

\subsection{Asymptotic Distribution : Dependent Random Variables}
Suppose that $\{X_{n}\}_{n\geq 1}$ is a sequence of random variables, and for any $a\in\mathbb{N}$ and $b\in\mathbb{N}$ such that $1\leq a\leq b <\infty$, let us denote ${\cal{F}}_{a}^{b} = \sigma(X_{a}, X_{a + 1}, \ldots, X_{b})$ and ${\cal{F}}_{m}^{\infty} = \sigma(X_{m}, X_{m + 1}, \ldots)$, where $\sigma(.)$ denotes the smallest $\sigma$-field generated by the random variables mentioned inside $(.).$ Note that if ${\cal{F}}_{1}^{k}$ and ${\cal{F}}_{k + n}^{\infty}$ are independent, then $P(A\cap B) - P(A) P (B) = 0$ for any $A\in {\cal{F}}_{1}^{k}$ and $B\in {\cal{F}}_{k + n}^{\infty}$. Let us quantify the dependence structure of the stationary sequence of random variables $\{X\}_{n\geq 1}$ using ${\cal{F}}_{1}^{k}$ and ${\cal{F}}_{k + n}^{\infty}$. Suppose that 
$$\psi(n) = \sup\left\{\left|\frac{P(A\cap B)}{P(A) P(B)} - 1\right| : A\in {\cal{F}}_{1}^{k}, B\in {\cal{F}}_{k + n}^{\infty}, k = 1, 2, \ldots\right\}.$$ Note that $\psi_{n} = 0$ for all $n\in\mathbb{N}$ if the sequence of random variables $\{X\}_{n\geq 1}$ are mutually independent, and the order of $\psi_{n}$ with respect to $n$ indicates the strength of the dependence among the random variables in the sequence $\{X_{n}\}_{n\geq 1}$. To understand the result for the dependent case, the following notations are introduced. For $c = 0, 1, \ldots, m$, let us denote 
$$\Phi_{c}(x_{1}, \ldots, x_{c}) = E\left\{\left(\frac{X_{1} + \ldots + X_{m}}{m}\right)\otimes\left(\frac{X_{1} + \ldots + X_{m}}{m}\right)|X_{1}, \ldots, X_{c}\right\},$$ and 
\begin{eqnarray}\label{GC}
g_{c}(x_{1}, \ldots, x_{c}) = \sum\limits_{d = 0}^{c} (-1)^{c - d}\sum\limits_{1\leq j_{1} <\ldots < j_{d}\leq c} \Phi_{d}(x_{j_{1}}, \ldots, x_{j_{d}}).
\end{eqnarray} In particular, note that $\Phi_{0} = E\left(\frac{X_{1} + \ldots + X_{m}}{m}\right)\otimes \left(\frac{X_{1} + \ldots + X_{m}}{m}\right)$ and $\Phi_{m} = \left(\frac{X_{1} + \ldots + X_{m}}{m}\right)\otimes\left(\frac{X_{1} + \ldots + X_{m}}{m}\right)$. Theorem \ref{depenasymp} describes the asymptotic distribution of $C_{m, n}$ for a certain dependent data, and the required conditions are the following. 

\vspace{0.1in}

\noindent (D1) $\{X_{n}\}_{n\geq 1}$ is a sequence of ${\cal{H}}$-valued identically distributed random elements, where ${\cal{H}}$ denotes the Hilbert space.

\vspace{0.1in}

\noindent (D2) $E||X_{1}||^{2}_{{\cal{H}}} < \infty$.

\vspace{0.1in}

\noindent (D3) $\sum\limits_{k = 1}^{\infty} (k + 1)^{m - 1} \{\psi(k)\}^{\frac{1}{2}} < \infty$.

\vspace{0.1in}

\noindent (D4) $\sigma^{2}_{\infty} := E ||g_{1}||^{2}_{{\cal{H}}\otimes  {\cal{H}}} + 2\sum\limits_{j = 2}^{\infty} E \langle g_{1}(X_{1}), g_{1}(X_{j}) \rangle_{{\cal{H}}\otimes {\cal{H}}}\neq 0$. 

\begin{theorem}\label{depenasymp} 
  Under (D1), (D2), (D3) and (D4),  $\frac{\sqrt{n}}{m}(C_{m, n} - C_{m})$ converges weakly to a ${\cal{H}}$-valued Gaussian random element with zero mean and the covariance operator $S$, where for $x\in {\cal{H}}$ and $y\in{\cal{H}}$, 
  \vspace{-0.1in}
\begin{eqnarray*}
  \langle Sx, y \rangle_{{\cal{H}}\otimes {\cal{H}}}
  &=& E \langle x, g_{1}(X_1) \rangle_{{\cal{H}}\otimes {\cal{H}}} \langle y, g_{1}(X_1) \rangle_{{\cal{H}}\otimes {\cal{H}}}
  + \sum\limits_{n = 2}^{\infty} E \langle x, g_{1}(X_1) \rangle_{{\cal{H}}\otimes {\cal{H}}}  \langle y, g_{1}(X_n) \rangle_{{\cal{H}}\otimes {\cal{H}}}\\
  &+& \sum\limits_{n = 2}^{\infty} E \langle x, g_{1}(X_n) \rangle_{{\cal{H}}\otimes {\cal{H}}}  \langle y, g_{1}(X_1) \rangle_{{\cal{H}}\otimes {\cal{H}}}.
  \end{eqnarray*} Here $g_{1}(.)$ is the same as $g_{c}(.)$ with $c = 1$ defined in \eqref{GC}.
\end{theorem}

\begin{remark}\label{redepasymp}
Condition (D1) is a minimal restriction for most of the Statistical methodologies, and this condition is easily verifiable for $\mathcal{H}$-valued random elements. For dependent sequence of random variables, (D3) indicates that the asymptotic distribution of $C_{m, n}$ depends on the order of $\psi(.)$, which measures the dependence among the $\mathcal{H}$-valued random elements. This condition is satisfied when $\psi (k) = O(\frac{1}{k^{2(l + 1)}})$, where $l\geq m$, and it implies that (D3) will be fulfilled for a wide range of dependence structure among the random elements.  The condition (D2) is mainly restriction on the moment on the normed random element, and it is also satisfied for many well-known random measures defined in Hilbert space such as Gaussian measure in Hilbert space. Condition (D4) implies that the random element $X$ is non-degenerate, and it leads to the asymptotic distribution as the Gaussian distribution. 
\end{remark}

\section{Numerical Study}\label{NS} As we mentioned in Section \ref{Intro}, here we carry out a small numerical study in a simple case to understand the effect of the choices of $m$. Let us consider ${\cal{B}} = \mathbb{R}$, and suppose that the data $X_{1}, \ldots, X_{n}$ are generated from standard $t$ distribution with $k$ degrees of freedom, where $k = 3, 4, \ldots, 10$. Here we are not considering $k = 1$ and 2 as the variance is not finite for such cases. Note that here $\theta(m) = 0$ as the standard $t$ distribution with $k$ ($3\leq k\leq 10$) degrees of freedom has mean $= 0$, and since ${\cal{B}} = \mathbb{R}$, we have 
$$C_{m, n} = \frac{1}{{n\choose m}}\sum_{1\leq i_{1}<i_{2}<\ldots<i_{m}\leq n} \left(\frac{X_{i_{1}} + \ldots + X_{i_{m}}}{m} - \theta (m)\right)^{2}$$ and 
$$C_{m} = E\left(\frac{X_{1} + \ldots + X_{m}}{m} - \theta (m)\right)^{2}.$$

In the numerical study, we consider $n = 10$, and $m = 1, \ldots, 10$. We generate data $L$ times with size $n$ and compute the empirical variance $EV(C_{m, n}) = \frac{1}{L}\sum\limits_{l = 1}^{L} (C_{m, n, l} - C_{m})$, where $C_{m, n, l}$ is the value of $C_{m, n}$ for the $l$-the generated data.   The values of $EV(C_{m, n})$ for various choices of $m\in\{1, \ldots, 10\}$ and $k\in\{3, \ldots, 10\}$ are reported in Table \ref{tab1} when $L = 100$.

\begin{table}[h!]
	
	\begin{center}		
		
		\begin{tabular}{|c|c|c|c|c|c|c|c|c|}\hline
			 & $k = 3$ & $k = 4$ & $k = 5$ & $k = 6$ & $k = 7$ & $k = 8$ & $k = 9$ & $k = 10$\\ \hline 
			$m = 1$ & {$2.88$} & {$2.01$} & {$1.67$} & {$1.54$} & {$1.39$} & {$1.33$} &  {$1.29$} & {$1.24$} \\ \hline
			$m = 2$ & {$2.75$} & {$1.94$} & {$1.58$} & {$1.43$} & {$1.38$} & {$1.35$} &  {$1.34$} & {$1.31$}\\ \hline
			$m = 3$ & {$2.66$} & {$1.91$} & {$1.53$} & {$1.44$} & {$1.41$} & {$1.38$} &  {$1.37$} & {$1.35$} \\ \hline
			$m = 4$ & {$2.54$} & {$1.88$} & 
			{$1.49$} & {$1.45$} & {$1.44$} & {$1.42$} &  {$1.39$} & {$1.38$} \\ \hline
			$m = 5$ & {$2.49$} & {$1.87$} & 
			{$1.48$} & {$1.46$} & {$1.43$} & {$1.42$} &  {$1.40$} & {$1.39$}\\ \hline
			$m = 6$ & {$2.50$} & {$1.89$} & 
			{$1.50$} & {$1.48$} & {$1.47$} & {$1.45$} &  {$1.44$} & {$1.44$}\\ \hline
			$m = 7$ & {$2.51$} & {$1.91$} & 
			{$1.52$} & {$1.50$} & {$1.48$} & {$1.47$} &  {$1.46$} & {$1.45$}\\ \hline
			$m = 8$ & {$2.59$} & {$1.96$} & 
			{$1.57$} & {$1.55$} & {$1.53$} & {$1.52$} &  {$1.50$} & {$1.49$}\\ \hline
			$m = 9$ & {$2.62$} & {$1.97$} & 
			{$1.59$} & {$1.58$} & {$1.56$} & {$1.55$} &  {$1.49$} & {$1.48$}\\ \hline
                 $m = 10$ & {$2.71$} & {$2.00$} & 
			{$1.63$} & {$1.55$} & {$1.48$} & {$1.45$} &  {$1.43$} & {$1.39$}\\ \hline
			
		\end{tabular}
	\end{center}
	\caption{\it The values (approximated up to two decimals) of $EV (C_{m, n})$ based on $L = 100$ when the data obtained from standard $t$-distribution with $k$ degrees of freedom. Here $n = 10$, $m\in \{1, \ldots, 10\}$ and $k\in \{3, \ldots, 10\}$.}
	
	\label{tab1}
	
\end{table}

The reported values in Table \ref{tab1} indicates that for heavy tailed distributions, i.e., when $k = 3, 4$ or 5, $C_{m, n}$ performs well $m = 5$ and 6, which is equal or close to $\frac{n}{2}$ (recall here $n = 10$) whereas for light tailed distributions, i.e., when $k = 10, 9$ or 8, $C_{m, n}$ performs well when $m = 1, 2$ or 10. To summarize, on the one hand, for light tailed distributions, the best performance comes when $C_{m, n}$ is essentially usual notion of covariance (i.e., $m = 1$), and then gradually becomes worsen till $m = \frac{n}{2}$. On the other hand, for heavy tailed distributions, the best performance of $C_{m, n}$ appears when $m$ is close to $\frac{n}{2}$. This observation may enable us to choose appropriate $m$, which controls the performance of $C_{m, n}$, in practice when the law of the data is unknown to us.  

\section{Concluding Remarks}\label{CR} 
This article investigates various large sample properties of a new covariance operator based on the concept of $U$ statistic for Banach valued random elements. Also, a small numerical study has been conducted to understand the choice of the order of kernels involved in the proposed co-variance operator.  This is essentially the summary of the work done in this article.

In the expression of $C_{m, n}$ (see \eqref{Sampleversion}) note that one may modify $C_{m, n}$ as
$$C_{m, n}^{k} = \frac{1}{{n\choose m}}\sum_{1\leq i_{1}<i_{2}<\ldots<i_{m}\leq n} k\left(\frac{X_{i_{1}} + \ldots + X_{i_{m}}}{m} - \theta (m)\right)\otimes k\left(\frac{X_{i_{1}} + \ldots + X_{i_{m}}}{m} - \theta (m)\right),$$ where $k : \mathcal{B}\rightarrow\mathcal{B}$ is a certain kernel. Note that $C_{m, n}^{k}$ coincides with $C_{m, n}$ when $k (x) = x$ for all $x\in\mathcal{B}$, and $k(x) = sign (x)$ for all $x\in\mathcal{B}$ leads to the sign covariance operator, which was studied by \cite{BOENTE2019115}. In this work, we have studied $C_{m, n}$ instead of $C_{m, n}^{k}$ to avoid additional notation complexity given the fact that the theoretical arguments would remain the same under certain conditions on $k$. Moreover, particularly, for the sign covariance operator, one can obtain the similar results for the sign covariance operator following the same arguments, and the derivation will be easier to a certain extent as $sign (x) = \frac{x}{||x||_{\cal{B}}}$ is a bounded function.

Recently there have been a few attempts to check whether two infinite dimensional random elements are independent or not (see, e.g., \cite{bhar2023testing} and a few references therein). One may be interested in investigating the same hypothesis problem using the proposed covariance operator in this article. Besides, as mentioned in Section \ref{Intro}, the proposed covariance operator can be used in functional principal component or outlier detection in a data.

We would like to close this section with the following discussion. In this work, we assume that the location of the data $X_{1}, \ldots, X_{n}$, i.e., $\theta(m)$ is known, and the results are derived using this fact. However, in a given problem, $\theta(m)$ may not be known beforehand, and to overcome it, one may replace $\theta(m)$ by its appropriate estimator. Deriving the similar results of the proposed covariance operator with unknown location may need more technicalities, and we will leave it as a future work.




\section{Appendix : Technical Details}\label{AP}


\noindent {\bf Proof of Proposition \ref{propu}:} 
Without loss of generality, we take $\theta (m) = 0$. From \eqref{Sampleversion},
\[C_{m, n} = \frac{1}{{n\choose m}}\sum_{1\leq i_{1}<i_{2}<\ldots<i_{m}\leq n} \left(\frac{X_{i_{1}} + \ldots + X_{i_{m}}}{m}\right)\otimes\left(\frac{X_{i_{1}} + \ldots + X_{i_{m}}}{m}\right).\]
For every $1\leq i_{1}<i_{2}<\ldots<i_{m}\leq n$, using both projective and injective tensor products, observe that
\[\left\| \frac{X_{i_{1}} + \ldots + X_{i_{m}}}{m} \otimes \frac{X_{i_{1}} + \ldots + X_{i_{m}}}{m}\right\|_{\cal{B} \otimes \cal{B}} \leq \left\| \frac{X_{i_{1}} + \ldots + X_{i_{m}}}{m} \right\|_{\cal{B}}^2\leq \frac{1}{m^2} \left(\sum_{j = 1}^m \|X_{i_{j}}\|_{\cal{B}} \right)^2.\]
Using i.i.d nature of the $X_i$'s, we have the following integrability condition
\[E\left\| \frac{X_{i_{1}} + \ldots + X_{i_{m}}}{m} \otimes \frac{X_{i_{1}} + \ldots + X_{i_{m}}}{m}\right\|_{\cal{B} \otimes \cal{B}} \leq \frac{m}{m^2} \sum_{j = 1}^m E \|X_{i_{j}}\|_{\cal{B}}^2= E||X||^{2}_{\cal{B}} < \infty.\]
Moreover, $\frac{X_{i_{1}} + \ldots + X_{i_{m}}}{m} \otimes \frac{X_{i_{1}} + \ldots + X_{i_{m}}}{m}$ has the same distribution as $\frac{X_1 + \ldots + X_m}{m} \otimes \frac{X_1 + \ldots + X_m}{m}$. Hence, the result follows.
\hfill$\Box$ 
\begin{lemma}\label{C2}
 Under (C2), for any $m\in\{1, \ldots, n\}$ and any $r\in\{1, \ldots, m\}$,
 
 \noindent $E\left\{\left\{\left(\frac{X_{1} + \ldots + X_{r - 1} + X_{r} + \ldots + X_{m}}{m} - \theta (m)\right)\otimes\left(\frac{X_{1} + \ldots + X_{r - 1} + X_{r} + \ldots + X_{m}}{m} - \theta (m)\right)\right\}|X_{1}, \ldots, X_{r - 1}\right\}$ 
 
 \noindent is a degenerate random element, and 
 
 \noindent $E\left\{\left\{\left(\frac{X_{1} + \ldots + X_{r} + X_{r+1} + \ldots + X_{m}}{m} - \theta (m)\right)\otimes\left(\frac{X_{1} + \ldots + X_{r}+ X_{r + 1} + \ldots + X_{m}}{m} - \theta (m)\right)\right\}|X_{1}, \ldots, X_{r}\right\}$
 
 \noindent is a non-degenerate random element.
\end{lemma}

\noindent {\bf Proof of Lemma \ref{C2}:} 
By Jensen's inequality for conditional expectation for Banach valued random elements (see, \cite{MR0576407}) and (C2), we have
\begin{align*}
&\left\| E\left.\left\{\left\{\left(\frac{X_{1} + \ldots + X_{r - 1} + X_{r} + \ldots + X_{m}}{m} - \theta (m)\right)\right.\right.\right.\right.\\
&\left.\left.\left.\otimes\left(\frac{X_{1} + \ldots + X_{r - 1} + X_{r} + \ldots + X_{m}}{m} - \theta (m)\right)\right\}|X_{1}, \ldots, X_{r - 1}\right\} \right\|  \\
&\leq E\left[\left\|\left\{\left(\frac{X_{1} + \ldots + X_{r - 1} + X_{r} + \ldots + X_{m}}{m} - \theta (m)\right)\right.\right.\right.\\
&\left.\left.\left.\otimes\left(\frac{X_{1} + \ldots + X_{r - 1} + X_{r} + \ldots + X_{m}}{m} - \theta (m)\right)\right\}\right\||X_{1}, \ldots, X_{r - 1}\right]   \\
&\leq E\left[\left\|\frac{X_{1} + \ldots + X_{r - 1} + X_{r} + \ldots + X_{m}}{m} - \theta (m)\right\|^2|X_{1}, \ldots, X_{r - 1}\right]\\
&=0.
\end{align*}
The degeneracy in the first part of the statement follows.

Since positivity of the variance implies the non-degeneracy of a random variable, we have the second part of the statement.
\hfill$\Box$

\begin{lemma}\label{C3}
Under (C3), for any $r\in\{1, \ldots, m\}$, 
$$E\left|\left|\left[E\left\{\left\{\left(\frac{X_{1} \ldots + X_{m}}{m} - \theta (m)\right)\otimes\left(\frac{X_{1} + \ldots + X_{m}}{m} - \theta (m)\right)\right\}|X_{1}, \ldots, X_{r}\right\}\right]\right|\right|^{\frac{mq}{q(m - r) + r}}_{\cal{B}\otimes\cal{B}} < \infty.$$
\end{lemma}

\noindent {\bf Proof of Lemma \ref{C3}:} 
Without loss of generality, we take $\theta (m) = 0$. Now,
\begin{align*}
& \left\|E\left\{\left\{\left(\frac{X_{1} + \ldots + X_{m}}{m}\right)\otimes\left(\frac{X_{1} + \ldots + X_{m}}{m}\right)\right\}|X_{1}, \ldots, X_{r}\right\}\right\|_{\cal{B}\otimes\cal{B}}\\
&\leq E\left\{\left\{\left\|\left(\frac{X_{1} + \ldots + X_{m}}{m}\right)\otimes\left(\frac{X_{1} + \ldots + X_{m}}{m}\right)\right\}\right\|_{\cal{B}\otimes\cal{B}} |X_{1}, \ldots, X_{r}\right\}\\
&\leq E\left\{\left\|\frac{X_{1} + \ldots + X_{r} + X_{r+1} + \ldots + X_{m}}{m}\right\|^2_{\cal{B}} |X_{1}, \ldots, X_{r}\right\}.
\end{align*}
Since, $\frac{mq}{q(m - r) + r} > 1$, by Jensen's inequality,
\begin{align*}
& E\left|\left|\left[E\left\{\left\{\left(\frac{X_{1} + \ldots + X_{m}}{m} - \theta (m)\right)\otimes\left(\frac{X_{1} + \ldots + X_{m}}{m} - \theta (m)\right)\right\}|X_{1}, \ldots, X_{r}\right\}\right]\right|\right|^{\frac{mq}{q(m - r) + r}}_{\cal{B}\otimes\cal{B}}\\
&\leq E \left( E\left\{\left\|\frac{X_{1} + \ldots + X_{r} + X_{r+1} + \ldots + X_{m}}{m}\right\|^2_{\cal{B}} |X_{1}, \ldots, X_{r}\right\}\right)^{\frac{mq}{q(m - r) + r}}\\
&\leq E \left( E\left\{\left\|\frac{X_{1} + \ldots + X_{r} + X_{r+1} + \ldots + X_{m}}{m}\right\|^{^{\frac{2mq}{q(m - r) + r}}}_{\cal{B}} |X_{1}, \ldots, X_{c}\right\}\right)\\
&= E \left\|\frac{X_{1} + \ldots + X_{r} + X_{r+1} + \ldots + X_{m}}{m}\right\|^{\frac{2mq}{q(m - r) + r}}_{\cal{B}}\\
&\leq  E \left(\frac{1}{m}\sum_{j = 1}^m \|X_j\|_{\cal{B}}\right)^{\frac{2mq}{q(m - r) + r}}\\
&\leq \frac{1}{m}  \sum_{j = 1}^m E\|X_j\|_{\cal{B}}^{\frac{2mq}{q(m - r) + r}}\\
&= E\|X_1\|_{\cal{B}}^{\frac{2mq}{q(m - r) + r}}
\end{align*}
The result follows.
\hfill$\Box$ 

The next result follows similar to Lemma \ref{C2}. We skip the proof to avoid repetitive arguments.
\begin{lemma}\label{C2*}
Under (C2*), $$E\left\{\left\{\left(\frac{X_{1} + \ldots + X_{m}}{m} - \theta (m)\right)\otimes\left(\frac{X_{1} + \ldots + X_{m}}{m} - \theta (m)\right)\right\}|X_{1}\right\}$$ is a non-degenerate random element.
\end{lemma}

Proof of Lemma \ref{C3*} follows similar to Lemma \ref{C3}. We skip the proof to avoid repetitive arguments.

\begin{lemma}\label{C3*} Under (C3*), for any $r\in\{1, \ldots, m\}$, $$E\left|\left|\left[E\left\{\left\{\left(\frac{X_{1} + \ldots + X_{m}}{m} - \theta (m)\right)\otimes\left(\frac{X_{1} + \ldots + X_{m}}{m} - \theta (m)\right)\right\}|X_{1}, \ldots, X_{r}\right\}\right]\right|\right|^{\frac{2r}{2r - 1}}_{\cal{B}\otimes\cal{B}} < \infty.$$
\end{lemma}


\begin{lemma}\label{D2}
Under (D1) and (D2), $$E \left|\left|\left(\frac{X_{1} + \ldots + X_{m}}{m}\right)\otimes\left(\frac{X_{1} + \ldots + X_{m}}{m}\right)\right|\right|_{{\cal{H}}\otimes{\cal{H}}} < \infty.$$
\end{lemma}

\noindent {\bf Proof of Lemma \ref{D2}:} 
Arguing similar to Proposition \ref{propu}, we have 
\[E\left\| \frac{X_{1} + \ldots + X_{m}}{m} \otimes \frac{X_{1} + \ldots + X_{m}}{m}\right\|_{\cal{H} \otimes \cal{H}} \leq \frac{m}{m^2} \sum_{j = 1}^m E \|X_{j}\|_{\cal{H}}^2= E||X||^{2}_{\cal{H}} < \infty.\]
\hfill $\Box$

\begin{lemma}\label{reversemar}
Let $B_{n} = \sigma(\omega : C_{m, n} (\omega), C_{m, n + 1} (\omega), \ldots)$, where $n\geq m$, and $\sigma(.)$ denotes the smallest sigma field formed by the collection of random variables mentioned inside (.). Then, the stochastic sequence $(C_{m, n}, B_{n})$ constitutes a reverse regular martingale, where 
$C_{m, n} = E\left\{(\frac{X_{1} + \ldots + X_{m}}{m} - \theta(m))\otimes (\frac{X_{1} + \ldots + X_{m}}{m} - \theta(m)) | B_{n}\right\}$ and $\theta(m) = E\left(\frac{X_{1} + \ldots + X_{m}}{m}\right)$.
\end{lemma}

\noindent {\bf Proof of Lemma \ref{reversemar}:} The proof follows the similar arguments provided in the proof of Lemma 1.1.3 in \cite{Borovskikh1996}. \hfill$\Box$

\begin{lemma}\label{ratehoeffconver}
For any $q\geq 1$, suppose that $||E(C_{m, k, c})||_{{\cal{B}}\otimes{\cal{B}}} < \infty$, where $C_{m, k, c}$ is defined in the proof of Theorem \ref{Rate_of_Convergence}, and let $\{b_{n}\}_{n\geq 1}$ be is a non-decreasing sequence of positive numbers. Then $$P\left(\max_{n\leq k\leq N} b_{k} E ||C_{m, k, c}||_{{\cal{B}}\otimes{\cal{B}}} > t\right)\leq\frac{1}{t^{q}}\left\{b_{n} E||C_{m, n}||^{q}_{{\cal{B}}\otimes{\cal{B}}} + \sum\limits_{k = n + 1}^{N} (b_{k} - b_{k - 1}) E||C_{m, k, c}||^{q}_{{\cal{B}}\otimes{\cal{B}}}\right\}$$ for all $c\leq n\leq N$ and for all $t > 0$. 
\end{lemma}

\noindent {\bf Proof of Lemma \ref{ratehoeffconver}:} The proof follows from the similar arguments provided in the proof of Theorem 2.3.1 in \cite{Borovskikh1996}. \hfill$\Box$ 

\noindent{\bf Proof of Theorem \ref{L1}} : Let $B_{n}$ be the same as defined in the statement of Lemma \ref{reversemar}. Note that $B_{n}$ is a decreasing sequence of set, i.e., $B_{n + 1}\subset B_{n}$ for all $n = m, m + 1, \ldots$, and $C_{m, n}$ is $B_{n}$-measurable random element. Now, denote $B_{\infty} = \cap_{n = m}^{\infty} B_{n}$, and using the assertion in Theorem 3 in \cite{Chatterji_1968}, which can be applied on ${\cal{B}}\otimes{\cal{B}}$-valued reverse regular martingale for convergence, along with the fact of Lemma \ref{reversemar}, we have $$C_{m, n}\rightarrow E\left\{\left(\frac{X_{1} + \ldots + X_{m}}{m} - \theta(m)\right)\otimes\left(\frac{X_{1} + \ldots + X_{m}}{m} - \theta(m)\right) | B_{\infty}\right\}$$ as $n\rightarrow\infty$ almost surely and in $L_{1}$ sense. Now, since $X_{1}, \ldots, X_{n}$ are i.i.d. sequence of random variables, in view of Hewitt-Savage theorem (see \cite{hewitt1955symmetric}), $B_{\infty}$ becomes trivial, i.e., $B_{\infty} = \{\emptyset, {\cal{X}}\}$, and hence, $$E\left\{\left(\frac{X_{1} + \ldots + X_{m}}{m} - \theta(m)\right)\otimes\left(\frac{X_{1} + \ldots + X_{m}}{m} - \theta(m)\right) | B_{\infty}\right\} = C_{m}.$$ It completes the proof. \hfill$\Box$

\noindent {\bf Proof of Corollary \ref{Consistency}}: For every $\epsilon > 0$, using Markov's inequality (see, e.g., \cite{van1998}), we have $$P[||C_{m, n} - C_{m}||_{{\cal{B}}\otimes {\cal{B}}} > \epsilon] < \frac{E[||C_{m, n} - C_{m}||_{{\cal{B}}\otimes {\cal{B}}}]}{\epsilon}.$$ Finally, the application on the assertion of Theorem \ref{L1} on the right hand side of the aforementioned inequality proves the result. \hfill$\Box$

\noindent {\bf Proof of Theorem \ref{Rate_of_Convergence}}: For any $\epsilon > 0$ and $1\leq q <2$, let us first define 
\begin{eqnarray}\label{relation1}
P_{n}(\epsilon, q) = P\left\{\omega :\sup_{k\geq n} {k\choose r}^{-\frac{1}{(q + 1)}}||C_{m, k}(\omega) - C_{m}||_{{\cal{B}}\otimes{\cal{B}}} \geq\epsilon \right\}.  \end{eqnarray} Now, using Hoeffding representation (see (1.1.9) in \cite{Borovskikh1996}), we have 
$$P_{n}(\epsilon, q) = P\left\{\omega :\sup_{k\geq n} {k\choose r}^{-\frac{1}{(q + 1)}}\left|\left|\sum\limits_{c = r}^{m} {n\choose c} C_{m, k, c}(\omega)\right|\right|_{{\cal{B}}\otimes{\cal{B}}} \geq\epsilon \right\},$$ where $C_{m, k, c} = \frac{1}{{n \choose c}}\sum\limits_{1\leq i_{1} < \ldots < i_{c}\leq n} E (C_{m, k}|X_{i_1}, \ldots, X_{i_c})$. Hence, using $P(\displaystyle\cup_{i = 1}^{r} A_{i})\leq\sum\limits_{i = 1}^{r} P (A_{i})$ for any arbitrary r many appropriate events $A_{1}, \ldots, A_{r}$, we have 
\begin{equation}\label{maininequality}
P_{n}(\epsilon, q)\leq\sum\limits_{c = r}^{m} P\left\{\omega: \sup_{k\geq n} {k\choose r}^{-\frac{1}{(q + 1)}}\left|\left| C_{m, k, c}(\omega)\right|\right|_{{\cal{B}}\otimes{\cal{B}}}\geq\epsilon_{c}\right\}, 
\end{equation} where $\epsilon_{c} = \frac{\epsilon}{(m - c + 1){m\choose c}}$. Let $2^{i -1}\leq n < 2^{i}$, and observe that 
\begin{eqnarray}\label{inequ1}
&&P\left\{\omega: \sup_{k\geq n} {k\choose r}^{-\frac{1}{(q + 1)}}\left|\left| C_{m, k, c}(\omega)\right|\right|_{{\cal{B}}\otimes{\cal{B}}}\geq\epsilon_{c}\right\}\nonumber\\
& = & P\left\{\omega: \bigcup_{k = n}^{\infty} {k\choose r}^{-\frac{1}{(q + 1)}}\left|\left| C_{m, k, c}(\omega)\right|\right|_{{\cal{B}}\otimes{\cal{B}}}\geq\epsilon_{c}\right\}\nonumber\\
 & \leq & P\left\{\omega: \bigcup_{j = i}^{\infty} \max_{2^{j - 1}\leq k\leq 2^{j}}{k\choose r}^{-\frac{1}{(q + 1)}}\left|\left| C_{m, k, c}(\omega)\right|\right|_{{\cal{B}}\otimes{\cal{B}}}\geq\epsilon_{c}\right\}\nonumber\\
 &\leq & \sum\limits_{j = i}^{\infty} P\left\{\omega: \max_{2^{j - 1}\leq k\leq 2^{j}}{k\choose r}^{-\frac{1}{(q + 1)}}\left|\left| C_{m, k, c}(\omega; \phi_{1})\right|\right|_{{\cal{B}}\otimes{\cal{B}}}\geq\frac{\epsilon_{c}}{2}\right\}\nonumber\\
 & + & \sum\limits_{j = i}^{\infty} P\left\{\omega: \max_{2^{j - 1}\leq k\leq 2^{j}}{k\choose r}^{-\frac{1}{(q + 1)}}\left|\left| C_{m, k, c}(\omega; \phi_{2})\right|\right|_{{\cal{B}}\otimes{\cal{B}}}\geq\frac{\epsilon_{c}}{2}\right\}, 
\end{eqnarray} where for $s = 1$ and 2, 

\begin{eqnarray}\label{ineq1*}
\phi_{s}(x_{1}, \ldots, x_{c}) = \sum\limits_{d = 1}^{c} (-1)^{c - d}\sum\limits_{1\leq j_{1}<j_{2}<\ldots < j_{d}\leq c}\phi_{s, c, d}(x_{j_1}, \ldots, x_{j_{d}}).
\end{eqnarray} Here $$\phi_{s, c, d}(x_{1}, \ldots, x_{d}) = C_{m, k, c}^{s} (x_{1}, \ldots, x_{d}) - E (C_{m, k, c}^{s} (x_{1}, \ldots, x_{d})),$$ 
where for any $\delta(c)\in (0, c]$, 
$$C_{m, k, c}^{1}(x_{1}, \ldots, x_{d}) = E\left(C_{m, k, c}(X_{1}, \ldots, X_{c})
1_{\{||C_{m, k, c}(X_{1}, \ldots, X_{c})||_{{\cal{B}}\otimes{\cal{B}}} \leq n^{\delta(c)}\}}|X_{1} = x_{1}, \ldots, X_{d} = x_{d}\right)$$ and $$C_{m, k, c}^{2}(x_{1}, \ldots, x_{d}) = E\left(C_{m, k, c}(X_{1}, \ldots, X_{c})
1_{\{||C_{m, k, c}(X_{1}, \ldots, X_{c})||_{{\cal{B}}\otimes{\cal{B}}} > n^{\delta(c)}\}}|X_{1} = x_{1}, \ldots, X_{d} = x_{d}\right).$$ 

Now, using the assertion in Lemma \ref{ratehoeffconver}, we have 
\begin{eqnarray}\label{ineq2}
&&P\left\{\omega: \max_{2^{j - 1}\leq k\leq 2^{j}}{k\choose r}^{-\frac{1}{(q + 1)}}\left|\left| C_{m, k, c}(\omega; \phi_{1})\right|\right|_{{\cal{B}}\otimes{\cal{B}}}\geq\frac{\epsilon_{c}}{2}\right\}\nonumber\\
&\leq & \frac{2^{p}}{\epsilon_{c}^{p}} \left({2^{j - 1}\choose r}^{\frac{p(q - 1)}{q}} E||C_{m, 2^{j - 1}, c}(\omega; \phi_{1})||^{p}_{{\cal{B}}\otimes{\cal{B}}}\right)\nonumber\\
& + & \sum\limits_{k = 2^{j - 1} + 1}^{2^{j}}\left({k\choose r}^{\frac{p(q - 1)}{q}} - {k - 1\choose r}^{\frac{p(q - 1)}{q}}  \right) E||C_{m, k, c}(\omega; \phi_{1})||_{{\cal{B}}\otimes{\cal{B}}}\nonumber\\
&\leq & a_{q}(c, r) 2^{-jc\left(\frac{p - \gamma_{rc}}{\gamma_{rc}}\right)}E \left(||C_{m, k, c}||^{p}_{{\cal{B}}\otimes{\cal{B}}}1_{(||C_{m, k, c}||_{{\cal{B}}\otimes{\cal{B}}} \leq 2^{jc}{\gamma_{rc}})}\right), 
\end{eqnarray} where $a_{p}(c, r)$ is a constant independent of $j$, and $\gamma_{rc}\in (1, p)$ is also a constant. Similarly, we also have 
\begin{eqnarray}\label{ineq3}
&&P\left\{\omega: \max_{2^{j - 1}\leq k\leq 2^{j}}{k\choose r}^{-\frac{1}{(q + 1)}}\left|\left| C_{m, k, c}(\omega; \phi_{2})\right|\right|_{{\cal{B}}\otimes{\cal{B}}}\geq\frac{\epsilon_{c}}{2}\right\}\nonumber\\
&\leq & \frac{2^{p}}{\epsilon_{c}^{p}} \left({2^{j - 1}\choose r}^{\frac{p(q - 1)}{q}} E||C_{m, 2^{j - 1}, c}(\omega; \phi_{2})||^{p}_{{\cal{B}}\otimes{\cal{B}}}\right)\nonumber\\
& + & \sum\limits_{k = 2^{j - 1} + 1}^{2^{j}}\left({k\choose r}^{\frac{p(q - 1)}{q}} - {k - 1\choose r}^{\frac{p(q - 1)}{q}}  \right) E||C_{m, k, c}(\omega; \phi_{2})||_{{\cal{B}}\otimes{\cal{B}}}\nonumber\\
&\leq & b_{q}(c, r) 2^{-jc\left(\frac{p - \gamma_{rc}}{\gamma_{rc}}\right)}E\left(||C_{m, k, c}||^{p}_{{\cal{B}}\otimes{\cal{B}}}1_{(||C_{m, k, c}||_{{\cal{B}}\otimes{\cal{B}}} > 2^{jc}{\gamma_{rc}})}\right), 
\end{eqnarray} where $b_{p}(c, r)$ is a constant independent of $j$. 

Therefore, using \eqref{inequ1}, \eqref{ineq2} and \eqref{ineq3}, we have 

\begin{eqnarray}\label{ineq4}
&&P\left\{\omega: \bigcup_{k = n}^{\infty} {k\choose r}^{-\frac{1}{(q + 1)}}\left|\left| C_{m, k, c}(\omega)\right|\right|_{{\cal{B}}\otimes{\cal{B}}}\geq\epsilon_{c}\right\}\nonumber\\  
 &\leq & \sum\limits_{j = i}^{\infty} P\left\{\omega: \max_{2^{j - 1}\leq k\leq 2^{j}}{k\choose r}^{-\frac{1}{(q + 1)}}\left|\left| C_{m, k, c}(\omega; \phi_{1})\right|\right|_{{\cal{B}}\otimes{\cal{B}}}\geq\frac{\epsilon_{c}}{2}\right\}\nonumber\\
 & + & \sum\limits_{j = i}^{\infty} P\left\{\omega: \max_{2^{j - 1}\leq k\leq 2^{j}}{k\choose r}^{-\frac{1}{(q + 1)}}\left|\left| C_{m, k, c}(\omega; \phi_{2})\right|\right|_{{\cal{B}}\otimes{\cal{B}}}\geq\frac{\epsilon_{c}}{2}\right\}\nonumber\\
 &\leq & \sum\limits_{j = i}^{\infty} [\{a_{q}(c, r) 2^{-jc\left(\frac{p - \gamma_{rc}}{\gamma_{rc}}\right)}E (||C_{m, k, c}||^{p}_{{\cal{B}}\otimes{\cal{B}}}1_{(||C_{m, k, c}||_{{\cal{B}}\otimes{\cal{B}}} \leq 2^{jc}{\gamma_{rc}})})\}\nonumber\\
 & + & \{b_{q}(c, r) 2^{-jc\left(\frac{p - \gamma_{rc}}{\gamma_{rc}}\right)}E (||C_{m, k, c}||^{p}_{{\cal{B}}\otimes{\cal{B}}}1_{(||C_{m, k, c}||_{{\cal{B}}\otimes{\cal{B}}} > 2^{jc}{\gamma_{rc}})})\}]\nonumber\\
 & = & \sum\limits_{j = 1}^{\infty}1_{(j\geq i)}[\{a_{q}(c, r) 2^{-jc\left(\frac{p - \gamma_{rc}}{\gamma_{rc}}\right)}E (||C_{m, k, c}||^{p}_{{\cal{B}}\otimes{\cal{B}}}1_{(||C_{m, k, c}||_{{\cal{B}}\otimes{\cal{B}}} \leq 2^{jc}{\gamma_{rc}})})\}\nonumber\\
 & + & \{b_{q}(c, r) 2^{-jc\left(\frac{p - \gamma_{rc}}{\gamma_{rc}}\right)}E (||C_{m, k, c}||^{p}_{{\cal{B}}\otimes{\cal{B}}}1_{(||C_{m, k, c}||_{{\cal{B}}\otimes{\cal{B}}} > 2^{jc}{\gamma_{rc}})})\}]\nonumber\\
 \end{eqnarray}

Now, in view of $1 < \gamma_{rc} < p$, observe that 
\begin{eqnarray}\label{ineq5}
&&\sum\limits_{j = 1}^{\infty} 1_{(j\geq i)} 2^{-jc\left(\frac{p - \gamma_{rc}}{\gamma_{rc}}\right)}E \left(||C_{m, k, c}||^{p}_{{\cal{B}}\otimes{\cal{B}}}1_{(||C_{m, k, c}||_{{\cal{B}}\otimes{\cal{B}}} \leq 2^{jc}{\gamma_{rc}})}\right)\nonumber\\  
&\leq & 2^{-ic\left(\frac{p - \gamma_{rc}}{\gamma_{rc}}\right)}\left(1 - 2^{-c\frac{p - \gamma_{rc}}{\gamma_{rc}}}\right)^{-1} + \left(1 - 2^{-c\frac{p - \gamma_{rc}}{\gamma_{rc}}}\right)^{-1}\nonumber\\  
&\times & \sum\limits_{k = 1}^{i - 1} 2^{-\frac{(i - k)c (p - \gamma_{rc})}{\gamma_{rc}}} E \left(||C_{m, k, c}||^{\gamma_{rc}}_{{\cal{B}}\otimes{\cal{B}}}1_{(2^{c(k - 1)}\leq||C_{m, k, c}||^{\gamma_{rc}}_{{\cal{B}}\otimes{\cal{B}}} \leq 2^{ck})}\right)\nonumber\\
& + & (1 - 2^{-c\frac{p - \gamma_{rc}}{\gamma_{rc}}})\sum\limits_{k = i}^{\infty} E \left(||C_{m, k, c}||^{\gamma_{rc}}_{{\cal{B}}\otimes{\cal{B}}}1_{(2^{c(k - 1)}\leq||C_{m, k, c}||^{\gamma_{rc}}_{{\cal{B}}\otimes{\cal{B}}} \leq 2^{ck})}\right).
\end{eqnarray}

As $2^{i - 1}\leq n < 2^{i}$, we have $2^{-\frac{ic(p - \gamma_{rc})}{\gamma_{rc}}}\leq n^{-\frac{c(p - \gamma_{rc})}{\gamma_{rc}}}$, and hence, 

\begin{eqnarray}\label{ineq6}
&&\sum\limits_{k = 1}^{i - 1} 2^{-\frac{(i - k)c (p - \gamma_{rc})}{\gamma_{rc}}} E\left(||C_{m, k, c}||^{\gamma_{rc}}_{{\cal{B}}\otimes{\cal{B}}}1_{(2^{c(k - 1)}\leq||C_{m, k, c}||^{\gamma_{rc}}_{{\cal{B}}\otimes{\cal{B}}} \leq 2^{ck})}\right)\nonumber\\
&\leq & 2 n^{-c\frac{p - \gamma_{rc}}{\gamma_{rc}}}E\left(||C_{m, k, c}||^{\gamma_{rc}}_{{\cal{B}}\otimes{\cal{B}}}1_{(||C_{m, k, c}||^{\gamma_{rc}}_{{\cal{B}}\otimes{\cal{B}}} > 1)}\right)\nonumber\\ 
&+& E\left(||C_{m, k, c}||^{\gamma_{rc}}_{{\cal{B}}\otimes{\cal{B}}}1_{(||C_{m, k, c}||^{\gamma_{rc}}_{{\cal{B}}\otimes{\cal{B}}} > n^{2c\frac{\gamma_{rc} - 1}{\gamma_{rc}}})}\right)
\end{eqnarray}
and 
\begin{eqnarray}\label{ineq7}
\sum\limits_{k = i}^{\infty} E \left(||C_{m, k, c}||^{\gamma_{rc}}_{{\cal{B}}\otimes{\cal{B}}}1_{(2^{c(k - 1)}\leq||C_{m, k, c}||^{\gamma_{rc}}_{{\cal{B}}\otimes{\cal{B}}} \leq 2^{ck})}\right)\leq E \left(||C_{m, k, c}||^{\gamma_{rc}}_{{\cal{B}}\otimes{\cal{B}}}1_{(2^{c}||C_{m, k, c}||^{\gamma_{rc}}_{{\cal{B}}\otimes{\cal{B}}} > n^{c})}\right).
\end{eqnarray}

Hence, using \eqref{ineq7}, \eqref{ineq6} and \eqref{ineq5}, we have 
\begin{eqnarray}\label{ineq8}
&&\sum\limits_{j = 1}^{\infty} 1_{(j\geq i)} 2^{-jc\left(\frac{p - \gamma_{rc}}{\gamma_{rc}}\right)}E \left(||C_{m, k, c}||^{p}_{{\cal{B}}\otimes{\cal{B}}}1_{(||C_{m, k, c}||_{{\cal{B}}\otimes{\cal{B}}} \leq 2^{jc}{\gamma_{rc}})}\right)\nonumber\\    
&\leq & 2 (2^{c(1 - \frac{1}{\gamma_{rc}})} - 1)^{-1} E \left(||C_{m, k, c}||^{p}_{{\cal{B}}\otimes{\cal{B}}}1_{(||C_{m, k, c}||_{{\cal{B}}\otimes{\cal{B}}} > n^{c})}\right)
\end{eqnarray}

As in \eqref{ineq4}, $a_{p}(c, r)$ and $b_{p}(c, r)$ are independent of $j$, using \eqref{ineq8}, \eqref{ineq4} and \eqref{maininequality} along with the conditions (C1), (C2) and (C3), we have 
\begin{eqnarray}\label{ineq9}
P_{n} (\epsilon, q)\leq\sum\limits_{c = r}^{m} a(c) n^{-c\mu(c)} + \sum\limits_{c = r}^{m} b(c) E\left(||C_{m, k, c}||^{\gamma_{rc}}_{{\cal{B}}\otimes{\cal{B}}}1_{(2^{c}||C_{m, k, c}||_{{\cal{B}}\otimes{\cal{B}}} > n^{\frac{2c(\gamma_{rc} - 1)}{\gamma_{rc}}})}\right),  
\end{eqnarray} where $a(c)$ and $b(c)$ do not depend on $n$, and $\mu_{c} = \min\left(\frac{p - \gamma_{rc}}{\gamma_{rc}}, (\frac{p - \gamma_{rc}}{\gamma_{rc}})^{2}\right)$.

Finally, the proof follows from \eqref{relation1},  \eqref{ineq9} and (C3). \hfill$\Box$

\noindent {\bf Proof of Theorem \ref{asymptotic_normality}:} First note that using Hoeffding representation (see (1.1.9) in \cite{Borovskikh1996}), we have 
\begin{eqnarray}\label{eq1}
&&\frac{\sqrt{n}}{m} (C_{m, n} - C_{m}) \nonumber\\
& = & \frac{1}{\sqrt{n}}\sum\limits_{i = 1}^{n} C_{m, n, 1}(X_{i}) + \frac{\sqrt{n}(m - 1)}{n(n - 1)}\sum\limits_{1\leq i < j\leq n} C_{m, n, 2} (X_{i}, X_{j}) + R_{n}, 
\end{eqnarray} where $$R_{n} = \frac{1}{m}\sum\limits_{c = 3}^{m}\frac{m!}{(m - c)!}\prod_{k = 1}^{c - 1} \left(1 - \frac{k}{n}\right)^{-1}n^{-\frac{c}{\gamma_{c}}}\sum\limits_{1\leq i_{1} <\ldots < i_{c}\leq n} C_{m, k, c} (X_{i_{1}}, \ldots, X_{i_{c}}). $$

Observe that $R_{n}$ depends on $m$, and $R_{n} = 0$ if $m = 1$ and 2. Moreover, for $c\geq 3$, (C3$^*$) indicates that one needs to have $\frac{6}{5}$-th conditional (conditioning of $X_{1}, \ldots, X_{c}$) moment of $(X_{1}+\ldots+X_{m})$. Furthermore, note that for $c = 3, \ldots, m$, the choice of $p$ described in (C1) (i.e., $\frac{3}{4}\leq p\leq 2$), $p > \gamma_{c} = \frac{2c}{2c - 1}$ as $\gamma_{c}$ is a decreasing function of $c$, and $\gamma_3 = \frac{6}{5}$. Then, it follows from the assertion in Theorem 3.1.2 in \cite{Borovskikh1996} that 
\begin{eqnarray}\label{reminder}
R_{n}\rightarrow 0   
\end{eqnarray} in probability as $n\rightarrow\infty$.

Next, note that for $c = 2$, $\gamma_{c} : = \frac{2c}{2c - 1} = \frac{4}{3}$. Therefore, again using Theorem 3.1.2 in \cite{Borovskikh1996}, we have 
$$n^{-\frac{3}{2}}\sum\limits_{1\leq i< j\leq n} C_{m, n, 2} (X_{i}, X_{j})\rightarrow 0$$ in probability as $n\rightarrow\infty$. Since $(m - 1)$ is independent of $n$, once can conclude that $\frac{\sqrt{n}(m - 1)}{n(n - 1)}\sum\limits_{1\leq i < j\leq n} C_{m, n, 2} (X_{i}, X_{j})\rightarrow 0$ almost surely as $n\rightarrow\infty$ when $p > \frac{4}{3}$.

Now, we want to analyse the limit of $n^{-\frac{3}{2}}\sum\limits_{1\leq i< j\leq n} C_{m, k, 2} (X_{i}, X_{j})$ as $n\rightarrow\infty$ when $p = \frac{4}{3} = \gamma_{2}$, i.e., $c = 2$. As $\delta(c)\in (0, c]$ (see below \eqref{ineq1*}), $\delta(2) = \frac{1}{2}$ is chosen. In \eqref{ineq1*}, using $c = 2$, $\delta(2) = \frac{1}{2}$, we have 
\begin{eqnarray}\label{eq1}
\phi_{1}(x_{1}, x_{2}) &=& C_{m, n, 2}(x_{1}, x_{2})1_{(||C_{m, n, 2}(x_1, x_2)||_{{\cal{B}}\otimes{\cal{B}}}\leq\sqrt{n})}\nonumber\\&& -  \int C_{m, n, 2}(x_{1}, y)1_{(||C_{m, n, 2}(x_1, y)||_{{\cal{B}}\otimes{\cal{B}}}\leq\sqrt{n})} P(dy)\nonumber\\ 
&& - \int C_{m, n, 2}(x, x_{2})1_{(||C_{m, n, 2}(x, x_{2})||_{{\cal{B}}\otimes{\cal{B}}}\leq\sqrt{n})} P(dx)\nonumber\\ 
&& + \int\int C_{m, n, 2}(x, y)1_{(||C_{m, n, 2}(x, y)||_{{\cal{B}}\otimes{\cal{B}}}\leq\sqrt{n})} P(dx) P(dy)
\end{eqnarray}

and 

\begin{eqnarray}\label{eq2}
\phi_{2}(x_{1}, x_{2}) &=& C_{m, n, 2}(x_{1}, x_{2})1_{(||C_{m, n, 2}(x_1, x_2)||_{{\cal{B}}\otimes{\cal{B}}} > \sqrt{n})}\nonumber\\&& -  \int C_{m, n, 2}(x_{1}, y)1_{(||C_{m, n, 2}(x_1, y)||_{{\cal{B}}\otimes{\cal{B}}} >\sqrt{n})} P(dy)\nonumber\\ 
&& - \int C_{m, n, 2}(x, x_{2})1_{(||C_{m, n, 2}(x, x_{2})||_{{\cal{B}}\otimes{\cal{B}}} > \sqrt{n})} P(dx)\nonumber\\ 
&& + \int\int C_{m, n, 2}(x, y)1_{(||C_{m, n, 2}(x, y)||_{{\cal{B}}\otimes{\cal{B}}} > \sqrt{n})} P(dx) P(dy)
\end{eqnarray}

Next, note that $\frac{1}{{n\choose 2}}\sum\limits_{1\leq i < j\leq n} C_{m, n, 2} (X_{i}, X_{j})$ is also a U-statistic, and one can write 
\begin{eqnarray}\label{eq3}
\frac{1}{{n\choose 2}}\sum\limits_{1\leq i < j\leq n} C_{m, n, 2} (X_{i}, X_{j}) & = & \frac{1}{{n\choose 2}}\sum\limits_{1\leq i < j\leq n} \phi_{1}(X_{i}, X_{j}) + \frac{1}{{n\choose 2}}\sum\limits_{1\leq i < j\leq n} \phi_{2}(X_{i}, X_{j}). 
\end{eqnarray}

Let's first work on the first term in the sum of the right hand side in \eqref{eq3}. As it is known in \cite{MR0515433}, using (2.2.1) in \cite{Borovskikh1996}, we have 
\begin{eqnarray}\label{ineq10}
E\left|\left|\frac{1}{{n\choose 2}}\sum\limits_{1\leq i < j\leq n} \phi_{1}(X_{i}, X_{j})\right|\right|_{{\cal{B}}\otimes{\cal{B}}}^{2}\leq\frac{\alpha_{2}}{{n\choose 2}^{2}} n^{\frac{1}{2}}\sum\limits_{i = 2}^{n} E\left|\left|\sum\limits_{j = 1}^{i - 1}\phi_{1}(X_{i}, X_{j})\right|\right|_{{\cal{B}}\otimes{\cal{B}}}^{2}, 
\end{eqnarray} where $\alpha_2$ is a constant. 

Now, for a fixed $X_{i} = x$, note that $\sum\limits_{j = 1}^{i - 1} \phi_{1}(x, X_{j})$ is the sum of independent ${\cal{B}}\otimes{\cal{B}}$-valued random elements such that $E(\phi_{1}(x, X_{j})) = 0$ for all $j = 1, \ldots, n$, and $E||\phi_{1}(x, X_{j})||_{{\cal{B}}\otimes{\cal{B}}} < \infty$ for all $j = 1, \ldots, n$. Since ${\cal{B}}\otimes{\cal{B}}$ is a Banach space of type-2 using (C1*), we have 
\begin{eqnarray}\label{ineq11}
E\left|\left|\sum\limits_{j = 1}^{i - 1}\phi_{1}(X_{i}, X_{j})\right|\right|_{{\cal{B}}\otimes{\cal{B}}}^{2}\leq b(i - 1) E||\phi_{1}(X_1, X_2)||_{{\cal{B}}\otimes{\cal{B}}},    
\end{eqnarray} where $b$ is a constant depending on ${\cal{B}}\otimes{\cal{B}}$. 

Next, using \eqref{ineq11} inside \eqref{ineq10}, we obtain 
\begin{eqnarray}\label{ineq12}
E\left|\left|\frac{1}{{n\choose 2}}\sum\limits_{1\leq i < j\leq n} \phi_{1}(X_{i}, X_{j})\right|\right|_{{\cal{B}}\otimes{\cal{B}}}^{2}
\leq\frac{\alpha_{2}b\sqrt{n}}{{n\choose 2}} E||\phi_{1}(X_1, X_2)||_{{\cal{B}}\otimes{\cal{B}}}. 
\end{eqnarray}

Further, note that 
\begin{eqnarray}\label{ineq13}
&&E||\phi_{1}(X_1, X_2)||_{{\cal{B}}\otimes{\cal{B}}}\leq 4^{2} E(||C_{m, n, 2}(X_{1}, X_{2})||_{{\cal{B}}\otimes{\cal{B}}}^{2}1_{(||C_{m, n, 2}(X_1, X_2)||_{{\cal{B}}\otimes{\cal{B}}}\leq\sqrt{n})})\nonumber\\
&&\leq 4^{2} E(||C_{m, n, 2}(X_{1}, X_{2})||_{{\cal{B}}\otimes{\cal{B}}}^{\frac{4}{3}}E(||C_{m, n, 2}(X_{1}, X_{2})||_{{\cal{B}}\otimes{\cal{B}}}^{\frac{2}{3}}1_{(||C_{m, n, 2}(X_1, X_2)||_{{\cal{B}}\otimes{\cal{B}}}\leq\sqrt{n})})\nonumber\\
&&\leq 4^{2}n^{\frac{1}{3}}E||C_{m, n, 2}(X_{1}, X_{2})||_{{\cal{B}}\otimes{\cal{B}}}^{\frac{4}{3}}.
\end{eqnarray}

Therefore, for $n\geq 2$, using \eqref{ineq13} inside \eqref{ineq12}, we have 
\begin{eqnarray}\label{ineq14}
E\left|\left|\frac{\sqrt{n}}{{n\choose 2}}\sum\limits_{1\leq i < j\leq n} \phi_{1}(X_{i}, X_{j})\right|\right|_{{\cal{B}}\otimes{\cal{B}}}^{2}\leq\frac{\alpha_{2}b 2^{6}}{n^{\frac{1}{6}}} E ||C_{m, n, 2}||^{\frac{4}{3}}_{{\cal{B}}\otimes{\cal{B}}}.   
\end{eqnarray} 

Hence, 
\begin{eqnarray}\label{ineq14*}
E\left|\left|\frac{\sqrt{n}}{{n\choose 2}}\sum\limits_{1\leq i < j\leq n} \phi_{1}(X_{i}, X_{j})\right|\right|_{{\cal{B}}\otimes{\cal{B}}}\rightarrow 0    
\end{eqnarray} as $n\rightarrow\infty$. 

Now, let us consider the second term in the sum of the right hand side in \eqref{eq3}. It follows from the assertion of Corollary 2.2.1 in \cite{Borovskikh1996} with $p =\frac{4}{3}$ in that result, we have \begin{eqnarray}\label{ineq15}
E\left|\left|\frac{1}{{n\choose 2}}\sum\limits_{1\leq i < j\leq n} \phi_{2}(X_{i}, X_{j})\right|\right|^{\frac{4}{3}}_{{\cal{B}\otimes{\cal{B}}}}\leq\alpha_{\frac{4}{3}}^{2}\times\frac{1}{{n\choose 2}^{\frac{1}{3}}}E||\phi_{2}(X_1, X_2)||^{\frac{4}{3}}_{{\cal{B}}\otimes{\cal{B}}}, 
\end{eqnarray} where 
\begin{eqnarray}\label{ineq16}
E||\phi_{2}(X_1, X_2)||^{\frac{4}{3}}_{{\cal{B}}\otimes{\cal{B}}}\leq 2^{\frac{8}{3}} E\left(||C_{m, n, 2}(X_{1}, X_{2})||^{\frac{4}{3}}_{{\cal{B}}\otimes{\cal{B}}}1_{(||C_{m, n, 2}(X_1, X_2)||_{{\cal{B}}\otimes{\cal{B}}} > \sqrt{n})}\right), \end{eqnarray} which follows from \eqref{eq2}.

Hence, for all $n\geq 2$, using \eqref{ineq16} inside \eqref{ineq15}, we have 
\begin{eqnarray}\label{ineq17}
&&E\left|\left|\frac{\sqrt{n}}{{n\choose 2}}\sum\limits_{1\leq i < j\leq n} \phi_{2}(X_{i}, X_{j})\right|\right|^{\frac{4}{3}}_{{\cal{B}\otimes{\cal{B}}}}\nonumber\\
& \leq & \alpha_{\frac{4}{3}}^{2} 2^{4} E\left(||C_{m, n, 2}(X_{1}, X_{2})||^{\frac{4}{3}}_{{\cal{B}}\otimes{\cal{B}}}1_{(||C_{m, n, 2}(X_1, X_2)||_{{\cal{B}}\otimes{\cal{B}}} > \sqrt{n})}\right).  
\end{eqnarray} Hence, 
\begin{eqnarray}\label{ineq18}
E\left|\left|\frac{\sqrt{n}}{{n\choose 2}}\sum\limits_{1\leq i < j\leq n} \phi_{2}(X_{i}, X_{j})\right|\right|^{\frac{4}{3}}_{{\cal{B}\otimes{\cal{B}}}}\rightarrow 0    
\end{eqnarray} as $n\rightarrow\infty$.

Therefore, using \eqref{ineq18} and \eqref{ineq14*} on \eqref{eq3}, we have 
\begin{eqnarray}\label{ineq18}
E\left|\left|\frac{\sqrt{n}}{{n\choose 2}}\sum\limits_{1\leq i < j\leq n} C_{m, n, 2} (X_{i}, X_{j})\right|\right|_{{\cal{B}}\otimes{\cal{B}}}\rightarrow 0 
\end{eqnarray} as $n\rightarrow\infty$. Afterwards, it follows from Chebyshev's inequality (see \cite{van1998}), we have 
\begin{eqnarray}\label{ineq19}
\left|\left|\frac{\sqrt{n}}{{n\choose 2}}\sum\limits_{1\leq i < j\leq n} C_{m, n, 2} (X_{i}, X_{j})\right|\right|_{{\cal{B}}\otimes{\cal{B}}}\rightarrow 0 
\end{eqnarray} in probability as $n\rightarrow\infty$. 

Finally, since ${{\cal{B}}\otimes{\cal{B}}}$ is a Banach space of type-2 (see (C1)) along with the fact \eqref{ineq19} and \eqref{reminder} applied on \eqref{eq3}, it follows from the central limit theorem of Banach valued random elements (see, e.g., \cite{MR0576407}), the result is proved. \hfill$\Box$

\noindent{\bf Proof of Theorem \ref{degenerateasymptoic}:} First note that using Hoeffding decomposition (see (1.1.9) in \cite{Borovskikh1996}), one has (see the derivation of \eqref{maininequality}),  
\begin{eqnarray}\label{eq5}
C_{m, n} - C_{m} &=& \sum\limits_{c = r}^{m} {m\choose c} C_{m, n, c}\nonumber\\
& = & {m\choose r} C_{m, n, r} + \sum\limits_{c = r + 1}^{m} {m\choose c} C_{m, n, c}.
\end{eqnarray}

Let us denote $\gamma_{rc} = \frac{cq}{q(c - r) + r}$, and note that $\gamma_{rc} < 2$ for all $c = r + 1, \ldots, m$. Therefore, using (C3), for $q = 2$ and in view of Theorem 3.1.2 in \cite{Borovskikh1996}, we have 
\begin{eqnarray}\label{lim1}
n^{\frac{r}{2}}\sum\limits_{c = r + 1}^{m} {m\choose c}C_{m, n, c}\rightarrow 0   
\end{eqnarray} in probability as $n\rightarrow\infty$. 

Moreover, it follows from Theorem 4.2.5 of \cite{Borovskikh1996} that 
\begin{eqnarray}\label{lim2}
n^{\frac{r}{2}}{m\choose r} C_{m, n, r}\rightarrow J_{r}(g_{r})    
\end{eqnarray} weakly as $n\rightarrow\infty$, where $J_{r}(g_{r})$ is the same as defined in \eqref{Ito}. Finally, the result follows from \eqref{lim1} and \eqref{lim2} using Slutsky's theorem (see, e.g., \cite{van1998}). \hfill$\Box$

\noindent {\bf Proof of Theorem \ref{depenasymp}:} Similar to \eqref{eq1}, first we observe that 
\begin{eqnarray}\label{eq6}
\frac{\sqrt{n}}{m}(C_{m, n} - C_{m}) = \frac{1}{\sqrt{n}}\sum\limits_{i = 1}^{n} C_{m, n, 1} (X_{i}) + R_{n},    \end{eqnarray} where $R_{n} = \frac{\sqrt{n}}{m}\sum\limits_{c = 2}^{m}{m\choose c}\frac{1}{{n\choose c}}\sum\limits_{1\leq i_{1}<\ldots < i_{c}\leq n} C_{m, k, c} (X_{i_{1}}, \ldots, X_{i_{c}})$.

To prove this result, we first try to show that 
\begin{eqnarray}\label{lim3}
R_{n}\rightarrow {\bf 0}    
\end{eqnarray} in probability as $n\rightarrow\infty$, where ${\bf 0}\in {\cal{B}}\otimes{\cal{B}}$. In order to establish \eqref{lim3}, let ${e_{j}}_{j\geq 1}$ be an orthonormal basis functions in ${\cal{H}}$, then we have 
\begin{eqnarray}\label{eq7}
&&E\left|\left|\frac{1}{{n\choose c}}\sum\limits_{1\leq i_{1}<\ldots < i_{c}\leq n} C_{m, k, c} (X_{i_{1}}, \ldots, X_{i_{c}})\right|\right|^{2}_{{\cal{H}}\otimes{\cal{H}}}\nonumber\\
& = & \sum\limits_{j = 1}^{\infty}E\langle e_{j}, \frac{1}{{n\choose c}}\sum\limits_{1\leq i_{1}<\ldots < i_{c}\leq n} C_{m, k, c} (X_{i_{1}}, \ldots, X_{i_{c}})\rangle^{2}_{{\cal{H}}\otimes{\cal{H}}}\nonumber\\
& = & \sum\limits_{j = 1}^{\infty}E\left(\sum\limits_{1\leq i_{1}<\ldots < i_{c}\leq n} \frac{1}{{n\choose c}}\langle e_{j}, C_{m, k, c} (X_{i_{1}}, \ldots, X_{i_{c}})\rangle^{2}_{{\cal{H}}\otimes{\cal{H}}}\right). 
\end{eqnarray} Now, note that $\left(\sum\limits_{1\leq i_{1}<\ldots < i_{c}\leq n} \frac{1}{{n\choose c}}\langle e_{j}, C_{m, k, c} (X_{i_{1}}, \ldots, X_{i_{c}})\rangle_{{\cal{H}}\otimes{\cal{H}}}\right)$ is a real valued $U$-statistic, and let us denote 
\begin{eqnarray}\label{eq8}
&&\Phi_{c}(x_1, \ldots, x_c)\nonumber\\
& =& E \left(\frac{X_{1} + \ldots + X_{m}}{m} - \theta(m)\right)\otimes\left(\frac{X_{1} + \ldots + X_{m}}{m} - \theta(m)\right) | X_{1} = x_{1}, \ldots, X_{c} = x_{c},  
\end{eqnarray} and it follows from \cite{Serfling1980} that 
\begin{eqnarray}\label{ineq20}
&&E\left(\sum\limits_{1\leq i_{1}<\ldots < i_{c}\leq n} \frac{1}{{n\choose c}}\langle e_{j}, C_{m, k, c} (X_{i_{1}}, \ldots, X_{i_{c}})\rangle^{2}_{{\cal{H}}\otimes{\cal{H}}}\right) \nonumber\\
&\leq & \frac{C}{n_{c}}\int\ldots\int\langle e_{j}, \Phi_{c}(x_{1}, \ldots, x_{c})\rangle^{2}_{{\cal{B}}\otimes{\cal{B}}}dP(x_{1})\ldots dP(x_{c}), 
\end{eqnarray} where $C$ is a some constant. Now, using \eqref{ineq20} on \eqref{eq7}, we have 
\begin{eqnarray}\label{ineq21}
&&E\left|\left|\frac{1}{{n\choose c}}\sum\limits_{1\leq i_{1}<\ldots < i_{c}\leq n} C_{m, k, c} (X_{i_{1}}, \ldots, X_{i_{c}})\right|\right|^{2}_{{\cal{H}}\otimes{\cal{H}}}\nonumber\\
& \leq &\frac{C}{n_{c}}\int\ldots\int\langle e_{j}, \Phi_{c}(x_{1}, \ldots, x_{m})\rangle^{2}_{{\cal{B}}\otimes{\cal{B}}}dP(x_{1})\ldots dP(x_{m}).    
\end{eqnarray} Hence, using (D2), we have 
\begin{eqnarray}\label{lim4}
E\left|\left|\frac{1}{{n\choose c}}\sum\limits_{1\leq i_{1}<\ldots < i_{c}\leq n} C_{m, k, c} (X_{i_{1}}, \ldots, X_{i_{c}})\right|\right|^{2}_{{\cal{H}}\otimes{\cal{H}}}\rightarrow 0     
\end{eqnarray} as $n\rightarrow\infty$. Therefore, by Chebyshev's inequality, we have 
\begin{eqnarray}\label{lim4}
\left|\left|\frac{1}{{n\choose c}}\sum\limits_{1\leq i_{1}<\ldots < i_{c}\leq n} C_{m, k, c} (X_{i_{1}}, \ldots, X_{i_{c}})\right|\right|_{{\cal{H}}\otimes{\cal{H}}}\rightarrow 0 \end{eqnarray} in probability as $n\rightarrow\infty$, and hence, 
\eqref{lim3} holds.

Now, we are working on $\frac{1}{\sqrt{n}}\sum\limits_{i = 1}^{n}C_{m, n, 1}(X_{i})$ (see \eqref{eq6}) as we have already established \eqref{lim3}. Let us now denote 
$\sigma_{n}^{2} = E\left|\left|\frac{1}{\sqrt{n}}\sum\limits_{i = 1}^{n}C_{m, n, 1}(X_{i})\right|\right|^{2}_{{\cal{B}}\otimes{\cal{B}}}$. Using straightforward algebra, we have 
\begin{eqnarray}\label{eq8}
\sigma_{n}^{2} = E||C_{m, n, 1}(X_{1})||^{2}_{{\cal{B}}\otimes{\cal{B}}} + 2\sum\limits_{j = 2}^{n} \left(1 - \frac{j}{n}\right) E \langle C_{m, n, 1}(X_{1}), C_{m, n, 1}(X_{j})\rangle_{{\cal{B}}\otimes{\cal{B}}}.   
\end{eqnarray} Note that under (D2) and (D3),  $\sigma_{n}^{2}\rightarrow\sigma_{\infty}^{2}$ (see (D4)) as $n\rightarrow\infty$. Hence, it follows from \cite{MR0576407} that $\frac{1}{\sqrt{n}}\sum\limits_{i = 1}^{n}C_{m, n, 1}(X_{i})$ converges weakly to the ${\cal{H}}$-valued Gaussian random element described in the statement of the theorem. Finally, using Stutsky's theorem on this fact and \eqref{eq7}, the result is proved. \hfill$\Box$

\section{Appendix: Topics on Banach spaces}

On Hilbert spaces, there is a natural choice of an tensor-product operation, which appears due to the underlying inner-product (see Remark \ref{Cmn-Hilbert}). Unlike the case in a Hilbert space, there are more than one choice of a tensor-product operation on a Banach space. We recall two such operations, called the projective tensor-product (see subsection \ref{A:projective-tensor}) and the injective tensor-product (see subsection \ref{A:injective-tensor}), from \cite{RyanBk}. In subsection \ref{A:type2-space}, we recall basic facts on type $p$ Banach spaces and $p$ uniformly smooth Banach spaces. Main references for this subsection are \cite{WoyczynskiBk, Ledoux-TalagrandBk}. Finally, in Subsection \ref{A:conditional-expectation}, we recall the definition of conditional expectation for Banach valued random elements from \cite[Chapter 2]{MetivierBk}.

In what follows, $\clb_1$ and $\clb_2$ are two real separable Banach spaces with norms denoted by $\|\cdot\|_1$ and $\|\cdot\|_2$.

\subsection{Projective tensor-product norm $\pi(.)$}\label{A:projective-tensor}

Let $u \in \clb_1 \otimes \clb_2$ with a representation $u = \sum_{j = 1}^n x_j \otimes y_j$ where $x_j \in \clb_1, y_j \in \clb_2, \forall j$. For any norm $\|u\|_{\clb_1 \otimes \clb_2}$, it is natural to expect the inequality
\[\|u\|_{\clb_1 \otimes \clb_2} \leq \sum_{j = 1}^n \|x_j\|_1 \|y_j\|_2.\]
This observation leads to the projective norm $\pi$ on $\clb_1 \otimes \clb_2$ defined as follows. For any $u \in \clb_1 \otimes \clb_2$, take
\begin{equation}\label{prj-norm-defn}
\pi(u) := \inf \left\{ \sum_{j = 1}^n \|x_j\|_1 \|y_j\|_2 : u = \sum_{j = 1}^n x_j \otimes y_j, \text{ with } x_j \in \clb_1, y_j \in \clb_2 \right\}.
\end{equation}
Note that $\pi(x \otimes y) = \|x\|_1 \|y\|_2$ for all $x \in \clb_1, y \in \clb_2$. 

It is also possible to describe the projective norm through bounded bilinear mappings on $\clb_1 \times \clb_2$. Let $B_{\clb_1}$ and $B_{\clb_2}$ denote the unit balls in $\clb_1$ and $\clb_2$, respectively. Now, consider the Banach space $\mathcal{B}(\clb_1 \times \clb_2)$ of real valued bounded bilinear mappings on $\clb_1 \times \clb_2$ with the norm 
\[\|B\| := \sup\{|B(x, y)| : x \in B_{\clb_1}, y \in B_{\clb_2}\}, \forall B \in \mathcal{B}(\clb_1 \times \clb_2).\]
Any $u = \sum_{j = 1}^n x_j \otimes y_j \in \clb_1 \otimes \clb_2$ acts on $B \in \mathcal{B}(\clb_1 \times \clb_2)$ by
\[\langle u, B \rangle = \sum_{j = 1}^n B(x_j, y_j).\]
Now the projective norm $\pi$ has the following identification (see \cite[p. 23]{RyanBk})
\[\pi(u) = \sup \{|\langle u, B \rangle| : B \in \mathcal{B}(\clb_1 \times \clb_2), \|B\| \leq 1\}.\]

\subsection{Injective tensor-product norm $\epsilon(.)$}\label{A:injective-tensor}

Let $u \in \clb_1 \otimes \clb_2$ with a representation $u = \sum_{j = 1}^n x_j \otimes y_j$ where $x_j \in \clb_1, y_j \in \clb_2, \forall j$. We associate $u$ with a real valued bounded bilinear form $B_u$ on $\clb_1^\ast \times \clb_2^\ast$ given by
\[B_u(\phi, \psi): = \sum_{j = 1}^n \phi(x_j) \psi(y_j), \, \forall \phi \in \clb_1^\ast, \psi \in \clb_2^\ast.\]
Let $B_{\clb_1^\ast}$ and $B_{\clb_2^\ast}$ denote the unit balls in $\clb_1^\ast$ and $\clb_2^\ast$, respectively. We define the projective norm $\epsilon(u)$ of $u \in \clb_1 \otimes \clb_2$ as the norm of $B_u$. Thus,
\begin{equation}\label{inj-norm-defn}
\epsilon(u) := \sup \left\{\left| \sum_{j = 1}^n \phi(x_j) \psi(y_j) \right| : \phi \in B_{\clb_1^\ast}, \psi \in B_{\clb_2^\ast}\right\}.
\end{equation}
The injective norm also has the following representations (see \cite[p. 46]{RyanBk})
\begin{equation}
\begin{split}
\epsilon(u) &= \sup \left\{\left\| \sum_{j = 1}^n \phi(x_j) y_j \right\|_2 : \phi \in B_{\clb_1^\ast}\right\}\\
&= \sup \left\{\left\| \sum_{j = 1}^n \psi(y_j) x_j\right\|_1 : \psi \in B_{\clb_2^\ast}\right\}.
\end{split}
\end{equation}

\subsection{Type $p$ Banach spaces and $p$ uniformly smooth Banach spaces}\label{A:type2-space}

In this subsection, we recall definitions and some basic facts on Type $p$ Banach spaces from \cite{WoyczynskiBk}. Another reference on this topic is \cite{Ledoux-TalagrandBk}. In what follows, $\clb$ shall denote a real separable Banach space with norm $\|\cdot\|_\clb$.

\begin{definition}[Rademacher Type $p$ spaces {\cite[Definition 6.2.1]{WoyczynskiBk}}]
Fix $p \geq 1$. Let $\{R_n\}_{n = 1}^\infty$ be a sequence of symmetric $\pm 1$ valued i.i.d. random variables. We say that the Banach space $\clb$ is of Rademacher Type $p$ if
\[\sup_{n} \inf\left\{\alpha \geq 0 : \forall x_1, x_2, \cdots, x_n \in \clb, \left(E \left\|\sum_{j = 1}^n R_j x_j \right\|_{\clb}^p \right)^{\frac{1}{p}} \leq \alpha \left( \sum_{j = 1}^n \left\|x_j \right\|_{\clb}^p \right)^{\frac{1}{p}}\right\} < \infty.\]
\end{definition}

\begin{definition}[Stable Type \texorpdfstring{p}{p} spaces {\cite[Definition 6.5.1]{WoyczynskiBk}}]
Fix $0 < p \leq 2$. Let $\{\xi_n\}_{n = 1}^\infty$ be a sequence of i.i.d. stable random variables with common Characteristic function $E e^{it\xi_1} = e^{-|t|^p}, \forall t \in \mathbb{R}$. We say that the Banach space $\clb$ is of Stable Type $p$ if
\[\sup_{n} \inf\left\{\alpha \geq 0 : \forall x_1, x_2, \cdots, x_n \in \clb, \left(E \left\|\sum_{j = 1}^n \xi_j x_j \right\|_{\clb}^{\frac{p}{2}}\right)^{\frac{2}{p}} \leq \alpha \left( \sum_{j = 1}^n \left\|x_j \right\|_{\clb}^p \right)^{\frac{1}{p}}\right\} < \infty.\]
\end{definition}

\begin{proposition}[{\cite[Proposition 7.1.1]{WoyczynskiBk}}]
A Banach space $\clb$ is of Rademacher Type 2 if and only if it is of Stable Type 2.
\end{proposition}

Following the above proposition, we now refer to Rademacher Type 2 or Stable Type 2 Banach spaces as just `Type 2' Banach spaces.

\begin{example}
By \cite[Theorem 7.6.1]{WoyczynskiBk}, due to \cite{KwapienPaper}, any Hilbert space $\cal{H}$ is a Type 2 Banach space. Observe that $\cal{H} \otimes \cal{H}$ is a Hilbert space, and is therefore a Type 2 Banach space. This is an example of the type of Banach spaces that we consider in this article.
\end{example}

\begin{definition}[$p$ Uniformly Smooth Banach spaces {\cite[Definition 3.1.2]{WoyczynskiBk}}]
Fix $1 < p \leq 2$. Define the modulus of smoothness of the Banach space $\clb$ as
\[\rho_{\clb}(\tau) := \sup\left\{\left\|\frac{x + y}{2}\right\|_{\clb} + \left\|\frac{x - y}{2}\right\|_{\clb} - 1 : \|x\|_\clb = 1, \|y\|_\clb = \tau\right\}, \tau \geq 0.\]
The Banach space $\clb$ is called $p$ uniformly smooth if $\rho_\clb(\tau) \leq C \tau^p, \forall \tau \geq 0$ for some constant $C$.
\end{definition}

\begin{example}[{\cite[Example 3.1.2]{WoyczynskiBk}}]
The $L^q, q > 1$ are $p$ uniformly smooth, where $p = \min\{q, 2\}$.
\end{example}

\begin{example}[{\cite[Corollary 3.1.1]{WoyczynskiBk}}]
Hilbert spaces are $2$ uniformly smooth.
\end{example}

\begin{remark}
Chapter 7 of \cite{WoyczynskiBk} contains an extensive survey on the results involving Kolmogorov's Three Series Theorem, the Law of Large Numbers, Central Limit Theorem and the Law of Iterated Logarithm. The interested reader may see the details here.
\end{remark}

\subsection{Conditional Expectation for Banach valued random elements}\label{A:conditional-expectation}

For real valued integrable random variables, the existence/definition of conditional expectation with respect to a sub-$\sigma$-field is derived from the Radon-Nikodym Theorem. Not all Banach spaces have the so-called Radon-Nikodym property, that the analog of Radom-Nikodym Theorem holds in the space (see \cite[Chapter 5]{RyanBk}). If the Banach space in consideration has the Radon-Nikodym property, then one can define the conditional expectation of Banach valued random elements in the usual way, as an application of the Radon-Nikodym Theorem.

If the Radon-Nikodym property is not true for the Banach space in question, there is an alternative approach in defining conditional expectation (see \cite[p. 41, 8.3 Remarks]{MetivierBk}).

Given $a_1, a_2, \cdots, a_n \in \clb$ and sets $A_1, A_2, \cdots, A_n \in \cal{A}$, consider the step function $f: \cal{X} \to \clb$ defined by $f(\omega) := \sum_{j = 1}^n 1_{A_j}(\omega) a_j, \forall \omega$. Using the notion of conditional expectation for real valued random variables with respect to a sub-$\sigma$-field $\cal{G}$, we have the conditional expectation $E [f \mid \cal{G}]$ given by
\[E [f \mid \mathcal{G}] := \sum_{j = 1}^n E [1_{A_j} \mid \mathcal{G}]\, a_j.\]

Since the above step functions form a dense subspace in $\mathcal{L}^1(\cal{X}, \cal{A}, P)$, we have the following densely defined continuous linear mapping $E [\cdot \mid \mathcal{G}] : \mathcal{L}^1(\cal{X}, \cal{A}, P) \to \mathcal{L}^1(\cal{X}, \cal{G}, P)$. Since this mapping is bounded on the said dense subspace, it extends uniquely to the whole of $\mathcal{L}^1(\cal{X}, \cal{A}, P)$. This provides an alternative approach in defining the conditional expectation for Banach valued random elements.

\noindent {\bf Acknowledgement:} Both authors are thankful to Professor BV Rao for stimulating discussion on fundamentals of the notion of covariance in infinite dimensional space. Subhra Sankar Dhar is grateful to Professor Yuliya Mishura for providing a useful reading material on U-statistic in Banach space, which gives the authors many ideas to derive the results. The work on Remark 2.2 is carried out motivated by a question from the audience in Joint Conference on Statistics and Data Science in China (JCSDS 2023) while Subhra Sankar Dhar was delivering a talk on covariance operators of Banach valued random elements in a session over there. Finally, Subhra Sankar Dhar gratefully acknowledges his core research grant (CRG/2022/001489), Government of India and Suprio Bhar acknowledges the support from the SERB MATRICS grant (MTR/2021/000517), Government of India.

\bibliography{lit} 

\end{document}